\documentclass{article}
\usepackage[english]{babel}
\usepackage[utf8]{inputenc}
\usepackage[T1]{fontenc}

\usepackage{amsfonts}
\usepackage{amsmath}
\usepackage{amssymb}
\usepackage{amsthm}
\usepackage{mathtools}
\usepackage{bbm}
\usepackage{geometry}
\usepackage[hyphens]{url}
\usepackage{hyperref}
\usepackage{nomencl}
\usepackage{authblk}

\usepackage{etoolbox}
\usepackage{blindtext}
\usepackage{url}
\usepackage{cleveref}
\usepackage[dvipsnames]{xcolor}
\usepackage{graphicx}
\usepackage{float}
\usepackage{subcaption}
\usepackage{algorithm} 
\usepackage{algpseudocode}

\hypersetup{
    colorlinks=true,
    linkcolor=black,
    anchorcolor=black,
    citecolor=black,
    filecolor=black,      
    urlcolor=black,
}
\DeclareMathOperator*{\argmax}{arg\,max}

\DeclareMathOperator*{\arginf}{arg\,inf}

\providecommand{\keywords}[1]{\textbf{\textbf{Keywords  }} #1}

\begin{document}

\title{Connecting Stochastic Optimal Control and Reinforcement Learning}

\date{August 2023}

\setcounter{Maxaffil}{0}
\renewcommand\Affilfont{\itshape\small}

\author[1]{Jannes Quer}
\author[1, 2]{Enric Ribera Borrell}
\renewcommand\Authands{ and }

\affil[1]{Institute of Mathematics, Freie Universität Berlin, 14195 Berlin, Germany}
\affil[2]{Zuse Institute Berlin, 14195 Berlin, Germany}

\maketitle

\begin{abstract}
In this paper the connection between stochastic optimal control and reinforcement learning is investigated. Our main motivation is to apply importance sampling to sampling rare events which can be reformulated as an optimal control problem. By using a parameterised approach the optimal control problem becomes a stochastic optimization problem which still raises some open questions regarding how to tackle the scalability to high-dimensional problems and how to deal with the intrinsic metastability of the system. To explore new methods we link the optimal control problem to reinforcement learning since both share the same underlying framework, namely a Markov Decision Process (MDP). For the optimal control problem we show how the MDP can be formulated. In addition we discuss how the stochastic optimal control problem can be interpreted in the framework of reinforcement learning. At the end of the article we present the application of two different reinforcement learning algorithms to the optimal control problem and a comparison of the advantages and disadvantages of the two algorithms.
\end{abstract}


\section{Introduction}
Rare event sampling is an area of research with applications in many different fields, such as finance \cite{Fournie1999}, molecular dynamics \cite{Quer2018} and many more. Very often the reason for the occurrence of rare events is that the dynamical system of interest exhibits metastable behaviour. Metastability means that the underlying process remains in certain regions of the state space for a very long time and only rarely changes to another region. This change is particularly important for accurate sampling of rare events. The average waiting time for the occurrence of the rare event is orders of magnitude longer than the timescale of the process itself. This behaviour is typically observed for dynamical systems following Langevin dynamics and moving in a potential with multiple minima. Here the metastable regions correspond to local minima of the potential. The minima are separated by barriers and the transitions between these regions are of interest. In molecular dynamics, for example, these quantities of interest correspond to the macroscopic properties of the molecules under consideration. Furthermore, it can be shown that the time to cross the barriers scales exponentially with the height of the barrier \cite{Berglund2013}. In terms of sampling it is observed that the variance of the Monte Carlo estimators associated with these rare transitions is often large. One idea to improve these estimators is the application of importance sampling but other methods such as splitting methods have been proposed. In this article we are going to focus on importance sampling. For a detailed discussion of splitting methods see \cite{Cerou2019} and the references therein. \\

One of the main challenges in importance sampling is to find a good bias so that the reweighted expectation has a low variance. The theory shows that the bias that would lead to a zero variance estimator is related to the quantity one wants to sample; see, e.g., \cite{Lelievre2010, Lelievre2016}. Therefore, many different variational methods have been proposed to find a good bias \cite{Valsson2014}. For importance sampling applications in stochastic differential equations it is well known that the optimal bias is actually given as the solution of a Hamilton-Jacobi-Bellman (HJB) equation, a nonlinear partial differential equation \cite{Fleming1993}. Since the HJB equation is the main equation of stochastic optimal control the importance sampling problem can be interpreted as stochastic optimal control. A stochastic optimization approach to approximate the bias using a parametric representation of the control has been proposed in \cite{Hartmann2012}. In the optimization approach the weights of the parametric representation are minimised to find the best approximation of the control. Methods for solving the related Hamilton-Jacobi-Bellman equation using deep learning based strategies in high dimensions have been developed; see, e.g., \cite{Weinan2017, Nusken2021, Weinan2021, Zhou2021, Han2020, Martin2022}. The approximation of control functions with tensor trains has been presented in \cite{fackeldey2020approximative}. \\

Although many methods have been proposed to approach the sampling problem from an optimal control point of view the stochastic optimization formulation offers the possibility to make a connection with reinforcement learning. Reinforcement learning (RL) is one of the three basic machine learning paradigms and has shown impressive results in high-dimensional applications such as Go and others \cite{Silver2016, Silver2017}. The reinforcement learning literature is very rich and many interesting ideas such as model-free, data-driven methods and robust gradient estimation have been extensively explored.    
From a more abstract point of view optimal control and reinforcement learning are concerned with the development of methods for solving sequential decision problems \cite{Powell2021} and their connection has been explored to some extent in \cite{Wang2020}. A general formulation can be given as follows: an intelligent agent should take different actions in an environment to maximise a so-called cumulative reward associated with a predefined goal. Applications of this formulation can be, for instance, a cleaning robot moving in a complex space, playing various games \cite{Sutton2018, Silver2014} or portfolio optimization \cite{Fleming1993}. The environment in which the agent moves is typically given in the form of a Markov decision process (MDP). Solution methods are often motivated by dynamic programming. The main difference between classical dynamic programming methods and reinforcement learning algorithms is that the latter do not require knowledge of an exact mathematical model of the MDP. In addition, reinforcement learning targets large MDPs where exact methods become infeasible. 
We make the link between the two fields in two ways. One way is to use the optimal control problem and formulate it as a Markov decision process which is the underlying theoretical framework of reinforcement learning. The other way is to formulate the reinforcement learning problem as a stochastic optimization problem. By comparing the resulting optimization problem with the optimization problem derived in \cite{Hartmann2012} for the importance sampling stochastic optimal control problem, one can see that both agree. Although the link is known there are few papers that show how to explicitly establish this link and the different language in these fields makes it difficult to understand the overlap. For this reason we want to show this connection in more detail for our problem of interest namely the application of importance sampling. However, the underlying connection is much more general and can be used for other applications that use stochastic modelling. Moreover, making this connection has the advantage that ideas developed in one area can be transferred to the other. \\

The paper is structured as follows. In \Cref{sec: is soc problem} we set the stage of rare event simulation, present the importance sampling problem and state its stochastic optimal control formulation. \Cref{sec: intro RL} is dedicated to the introduction of the reinforcement learning framework. In \Crefrange{subsec: MDP}{subsec: trajectories return and value functions} we discuss Markov decision problems, which are the underlying theoretical framework of reinforcement learning. In \Cref{subsec: rl optimization problem} we state the reinforcement learning optimization problem and in \Cref{subsec: rl algorithms} we recap the key ideas behind the main RL algorithms. \Cref{sec: soc problem as rl} is devoted to showing how optimal control and reinforcement learning are related. In \Cref{subsec: is soc rl environment} we show how the stochastic optimal control problem can be formulated as an MDP. In \Cref{subsec: soc vs rl} we compare the two optimization formulations presented for both problems. In \Cref{subsec: algorithms for soc} we first discuss how a previously presented solution method in the framework of reinforcement learning can be understood as a model-based approach. We introduce the well-known model-free Deterministic Policy Gradient (DPG) algorithms. In \Cref{sec: numerical examples} we present an application of the presented methods to a small toy problem. The focus of this section is the approximation of the optimal control problem which we will discuss in more detail. We conclude the article with a summary of our results.

\section{Importance sampling SOC problem}
\label{sec: is soc problem}
The main goal of this paper is to show how stochastic optimal control (SOC) and reinforcement learning (RL) are related. We introduce a SOC problem related to an importance sampling diffusion problem. However, the importance sampling SOC formulation can be easily adapted to general SOC problems and is not restricted to the importance sampling application. First we motivate the importance sampling problem and present a rather formal formulation. Then we show the relationship to the corresponding Hamilton-Jacobi-Bellman (HJB) equation. Finally we show how this problem can be formulated as a stochastic optimization problem. This section will be formulated in continuous time to show the relation to the HJB equation. In the following sections we will return to the discrete-time problem for convenience as this is the viewpoint used in most of the reinforcement learning literature.

\subsection{The sampling problem}
The problem we are going to consider in this paper is a particle moving in some potential landscape $V: \mathbb{R}^d \rightarrow \mathbb{R}$ on the bounded domain $\mathcal{D} \subset \mathbb{R}^d$. The movement of the particle follows the gradient of the potential plus a perturbation given by a Brownian motion noted as $w_t$ and  scaled by a diffusion term. This type of movement can be described by a stochastic differential equation (SDE) known as the overdamped Langevin equation: 
\begin{equation}
\label{eq: sde}
\mathrm{d}x_t = -\nabla V(x_t) \mathrm{d}t + \sigma(x_t) \mathrm{d}w_t, \quad x_0 = x,
\end{equation}
where $x_t$ is the position of the particle at time $t$, $x \in \mathbb{R}^d$ is a deterministic starting position and the diffusion term is chosen to be $\sigma(x) = \sqrt{2 \beta^{-1}}$ \text{Id}. Let us assume we are interested in computing the expectation
\begin{equation}
\label{eq: psi}
\Psi(x) \coloneqq \mathbb{E} \bigl[I(x_{0:T}) \mid x_0=x \bigr]
\end{equation}
of path dependent characteristics which can be described as
\begin{equation}
\label{eq: quantity of interest}
I(x_{0:T}) \coloneqq \exp \Bigl(- g(x_T) - \int_0^T f(x_t) \mathrm{d}t \Bigr),
\end{equation}
where $T = \inf\{t>0 \mid x_t \in \mathcal{T} \}$ is the first hitting time of a specific target set $\mathcal{T} \subset \mathcal{D}$ and $f, g: \mathbb{R}^d \rightarrow  \mathbb{R}$ are some bounded and sufficiently smooth functions. In the literature $f$ and $g$ are often denoted as the \textit{running cost} and the \textit{terminal cost} respectively. For instance the quantity of interest which leads to estimating the moment generating function of $T$ can be achieved by setting $f=\lambda$ and $g=0$. 

We usually consider a potential $V$ which has many local minima and which induces metastable dynamics in the model under consideration. This metastability leads to a high variance of the Monte Carlo estimator of the path dependent quantity of interest and thus it becomes unreliable. To improve the properties of the estimator we will apply an importance sampling strategy. The general idea of such a strategy is to sample a different dynamical system by introducing a bias and later correcting this effect in the expectation. \\

In the considered problem we are going to introduce a bias in the drift term of the SDE. This additional drift is known as a control in the SOC literature. By construction it does not influence the stochastic noise. The controlled dynamical system is given by 
\begin{equation}
\label{eq: controlled sde}
\mathrm{d}y_t = \bigl(-\nabla V(y_t) + \sigma(y_t) u(y_t) \bigr) \mathrm{d}t + \sigma(y_t) \mathrm{d}w_t, \quad y_0 = x,
\end{equation}
where $u: \mathbb{R}^d \rightarrow \mathbb{R}^d$ is the mentioned control, which belongs to the space of time-independent feedback controls $\mathcal{U}$ (see \cite{Hartmann2017} for more technical assumptions on the control). The time-independence assumption on the feedback control comes from the fact that the quantity of interest \eqref{eq: psi} is time independent thus the overall problem we are trying to solve is stationary. The controlled dynamics can be related to the original dynamics by using the Girsanov formula, which is a change of measure in path space. Precisely, we have that the importance sampling quantity of interest is an unbiased estimator of \eqref{eq: psi}
\begin{equation}
\label{eq: is estimator}
\Psi(x) = \mathbb{E} \bigl[I(x_{0:T}) \mid x_0=x \bigr] = \mathbb{E} \bigl[I(y_{0:T_u}) \, m_{0:T_u} \mid y_0=x \bigr],
\end{equation}
where $m_{0:T_u}$ is an exponential Martingale given by
\begin{equation}
m_{0:T_u} \coloneqq \exp{\left(-  \int_0^{T_u} u(y_t) \cdot \mathrm{d}w_t - \frac{1}{2} \int_0^{T_u} |u(y_t)|^2 \mathrm{d}t \right)}
\end{equation}
and ${T_u}$ is the first hitting time of $\mathcal{T}$ under the controlled process. Although the importance sampling relation \eqref{eq: is estimator} holds for any control $u \in \mathcal{U}$, the variance of the corresponding estimators significantly depends on the choice of $u$. Hence, one is tempted to aim for an optimal control $u^*$ which minimizes the variance of the importance sampling estimator over the space of controls
\begin{equation}
\label{eq: variance minimization}
u^* = \arginf\limits_{u \in \mathcal{U}} \biggr\{ \mathrm{Var} \Bigr( I(y_{0:T_u}) \, m_{0:T_u} \mid y_0=x \Bigl) \biggl\}.
\end{equation}

\subsection{Hamilton-Jacobi-Bellman equation}
With this formulation at hand we can now search for the optimal control. It is well-known in the literature that the optimal control can be found by using the Feynman-Kac formula and the resulting partial differential equation (PDE) connection \cite{Lelievre2016, Hartmann2017}. Via the Feynman-Kac formula $\Psi$ satisfies the following elliptic boundary value problem (BVP)

\begin{subequations}
\label{eq: linear hjb}
\begin{align}
(\mathcal{L} -f) \Psi &= 0 \qquad \qquad \, \text{in} \; \mathcal{O} \\
\Psi &= \exp(-g) \quad \text{on} \; \partial{\mathcal{O}},
\end{align}
\end{subequations}
on the domain $\mathcal{O} \coloneqq \mathcal{D} \cap \mathcal{T}^c$ where $\mathcal{L}$ denotes the infinitesimal generator, which is given by
\begin{equation*}
\mathcal{L} = \frac{1}{2} \sum_{i,j=1}^d (\sigma \sigma^\top)_{ij}(x) \frac{\partial^2}{\partial x_i \partial x_j} - \sum_{i=1}^d \Bigl(\frac{\partial}{\partial x_i}V(x) \Bigr) \frac{\partial}{\partial x_i} .
\end{equation*}
By using the Cole-Hopf transformation $\Phi = - \log \Psi$, one can derive the well-known Hamilton-Jacobi-Bellman equation
\begin{subequations}
\label{eq: hjb}
\begin{align}
\mathcal{L}\Phi - \frac{1}{2} |\sigma^\intercal \nabla \Phi|^2 + f\Phi = 0 \quad &\text{in} \; \mathcal{O} \\
\Phi = g \quad &\text{on} \; \partial{\mathcal{O}}.
\end{align}
\end{subequations}
Moreover, it has been shown that the optimal control depends directly on the solution of the above PDE
\begin{equation}
\label{eq: optimal control characterization}
u^* 
= - \sigma^T \nabla \Phi 
= \sigma^T \nabla \log{\Psi}
\end{equation}
and that the corresponding importance sampling estimator achieves zero variance see e.g. \cite{Lelievre2016, Hartmann2017}. \\

A priori one can calculate the quantity that we originally wanted to estimate as a function of the initial position via e.g. a finite difference method. However, such a problem becomes not trivial to solve for high-dimensional settings due to the curse of dimensionality. Furthermore, we know that the problem we are trying to solve is hard because we are trying to find a solution for a nonlinear PDE.

\subsection{Stochastic optimization problem}
To this end we are going to reformulate the problem as an optimization problem. Let us recall the value function $\Phi$ as a function of $\Psi$ which can be expressed in terms of the importance sampling estimator \eqref{eq: is estimator}
\begin{equation*}
\Phi(x) 
= - \log{\left( \mathbb{E} \left[ \exp{\left(
- g(y_{T_u})
- \int_0^{T_u} f(y_t) \mathrm{d}t 
- \int_0^{T_u} u(y_t) \cdot \mathrm{d}w_t 
- \frac{1}{2} \int_0^{T_u} |u(y_t)|^2 \mathrm{d}t
\right)} \Bigm| y_0 = x \right] \right)}
\end{equation*}
By using the Jensen's inequality and later the fact that an expectation of an It\^o stochastic integral is zero we obtain the following upper bound for the value function
\begin{equation}
\label{eq: value function upper bound}
\Phi(x) 
\leq \mathbb{E} \left[
g(y_{T_u})
+ \int_0^{T_u} f(y_t) \mathrm{d}t 
+ \frac{1}{2} \int_0^{T_u} |u(y_t)|^2 \mathrm{d}t \Bigm | y_0 = x \right].
\end{equation}
In the literature it is shown with the help of some stochastic calculus that the inequality for the optimal control $u^*$ is obtained \cite{Hartmann2017, Lelievre2016}.
We can understand the right-hand side of the equation \eqref{eq: value function upper bound} as a performance measure of the control applied to the dynamical system. The first two terms in the expectation are equal to the logarithm of the quantity of interest e.g. the stopping time in our example setting and the third term measures the force applied to the system. In the optimal control literature this expression is known as the \textit{cost} functional conditional on the initial condition. On the one hand we can use the right-hand side of \eqref{eq: value function upper bound} as the objective function for the optimization problem since it is easier to handle numerically. On the other hand we see that there is still a connection to the optimal control problem. If we manage to find the optimal solution with the optimization procedure we have the solution to the HJB equation. \\

By minimizing this functional in terms of the control $u$ subject to the controlled dynamical system we have derived an optimization problem in order to find the optimal control.
This optimal control problem is given by 
\begin{equation}
\label{eq: soc problem}
u^* = \arginf\limits_{u \in \mathcal{U}} J(u; x)
\end{equation}
where the minimisation is taken again over the space of controls and $J(u; x)$ is the corresponding cost functional
\begin{subequations}
\begin{align} 
J(u; x) &\coloneqq \mathbb{E}\left[ g(y_{T_u}) + \int_0^{T_u} f(y_t) \mathrm{d}t + \frac{1}{2} \int_0^{T_u} |u(y_t)|^2 \mathrm{d}t \Bigm| y_0 = x \right] \label{eq: soc objective function} \\ 
\text{s.t.}\ \mathrm{d}y_t &=  \bigl(-\nabla V(y_t) + \sigma(y_t) u(y_t) \bigr) \mathrm{d}t + \sigma(y_t) \mathrm{d}w_t, \quad y_0 = x . \label{eq: soc dynamics}
\end{align}
\end{subequations}
Many different ideas have been proposed in the literature to solve this problem. One idea is to use a Galerkin projection of the control into a space of weighted finite initial functions and optimise over the weights using a gradient descent method or a cross-entropy method \cite{Hartmann2012, Hartmann2016}. Another idea is to solve the deterministic control problem and use it to steer the dynamical system in the right direction \cite{Weare2012}. A more detailed discussion of solution methods can be found in \cite{Lelievre2016}.

According to the theory using optimal control to sample the quantity of interest would result in a zero variance estimator. In principle one needs to sample a trajectory to find the quantity of interest. However, due to discretization and numerical issues this is not possible in the implementation. However, with a good approximation to the optimal control the sampling effort can be massively reduced and the estimator converges faster with a smaller relative error. Furthermore, it can be shown that the relative error scales exponentially with the approximation error \cite{Richter2021}. This means that one should use the best possible approximation to find a good estimator of the quantity of interest. Because of this dependence, it is necessary to use methods that find a solution close to the optimum. \\

For mathematical completeness it is necessary to have a small remark on other possible types of time horizons. In this article $T_u$ is an a.s. finite stopping time with respect to the canonical filtration of the controlled process. In general we could consider 
\begin{itemize}
\item a finite horizon time $T_u=T_\text{end}$, leading to a deterministic stopping time,
\item a bounded stopping times $\widetilde{T_u} = \min(T_u, T_\text{end})$,
\item or a general random stopping times $T_u \in \mathbb{R^+} \cup \{+\infty\}$.
\end{itemize}
The first two types of stopping times are also a.s. finite and are of special interest since they guarantee the applicability of Girsanov's theorem for a control on a bounded domain. In this case the Novikov condition is satisfied. The reason is that with this assumption the boundedness of the control is guaranteed.
However, one ends up with a different problem formulation. In these cases the corresponding BVPs are not elliptic but parabolic leading to time-dependent solutions and as a consequence the optimal control is also time-dependent \cite{Hartmann2014, Hartmann2017}. The general case of random stopping time needs to be discussed in more detail and is beyond the scope of this article \cite{Lelievre2016}.

\section{Introduction to reinforcement learning}
\label{sec: intro RL}
Before showing the connection between RL and SOC for completeness we would like to give a brief overview of a typical model for a reinforcement learning problem. First, we look at the underlying theoretical framework namely Markov decision processes. Then, we will discuss key concepts of RL theory, such as types of policies, value and Q-value functions, and their recursive relations provided by the Bellman equations. Finally, we will introduce the reinforcement learning problem from an optimization point of view and discuss two different formulations of the RL problem. \\

As stated by Sutton and Barto: ``Reinforcement learning is learning what to do - how to map situations to actions - in order to maximise a numerical reward signal.'' \cite{Sutton2018}. So from this very first definition we can already see what it takes to define a reinforcement learning problem. To learn what ``to do'', we need someone to do something. This is usually called the \textit{agent}. The agent experiences situations by interacting with an \textit{environment}. This interaction is based on the \textit{actions} the agent takes, the \textit{reward signals} the agent receives from the environment and the \textit{next state} the agent has moved to by following the controlled dynamics. Usually the interaction with the environment is considered to take place over some time interval (finite or infinite). \\

The goal of the agent is to reach a predefined goal that is part of the environment. The agent will reach the goal if it optimally chooses the sequence of actions so that the sum of the reward signals received along a trajectory is maximised. Given the dynamics of the environment and a chosen goal reinforcement learning assumes that there are reward functions that can lead the agent to success in the sense that the predefined goal is reached. The existence of a unique reward function that is certain to lead to the predefined goal is not given. Thus, the choice of the reward function can be flexible and depend on the task that the agent has to solve. Furthermore, it may have an influence on the learned action.

Environments with sparse reward functions are often difficult to handle, since the reward signal received is often negligible. For example this is often the case in games where the reward signal is $0$ throughout the game and changes to $-1$ in the event of a loss or $+1$ in the event of a win. However, this is not the only possible reward function that makes the agent learn how to play. If a person plays a game of chess with a teacher who tells her how she played after every move the person will receive a much richer reward signal which will eventually lead to faster learning.

\subsection{Markov decision processes}
\label{subsec: MDP}
The theoretical framework of all reinforcement learning problems is a Markov decision process. A typical RL problem has the following elements:
\begin{itemize}
\item the \textit{state space} $\mathcal{S}$ is the set of states i.e. the set of all possible situations in the environment
\item the \textit{action space} $\mathcal{A}$ is the set of actions the agent can choose at each state
\item the \textit{set of decision epochs} $\mathbb{T} \subset \mathbb{R}^+$ is the set of time steps corresponding to the times where the agent acts. Let us assume the set of decision epochs is discrete. In this case, it can either be finite, i.e., $\mathbb{T} = \{ 0, \dots, T\}$ with $T \in \mathbb{N}^+$, or infinite $\mathbb{T} = \mathbb{N}_0$
\item the \textit{(state-action) transition probability function} $p: \mathcal{S} \times \mathcal{S} \times \mathcal{A} \rightarrow [0, 1]$ provides the probability of transitioning between states after having chosen a certain action. The transition probability function depends on the state where the agent is $s \in \mathcal{S}$, the action she chooses $a \in \mathcal{A}$ and the next state the agent moves into $s' \in \mathcal{S}$, i.e. 
\begin{equation}
\label{eq: state transition probability}
p(s', s, a) = p(s' | s, a) \coloneqq \mathbb{P}(s_{t+1} = s' \mid s_t=s, \, a_t=a ).
\end{equation}
The transition probability function given a state-action pair $p(\cdot| s, a)$ is a conditional probability mass function
\begin{equation}
\sum\limits_{s' \in \mathcal{S}} p(s'|s, a) = 1.
\end{equation}
If $\mathcal{S}$ is a continuous state space $p: \mathcal{S} \times \mathcal{S} \times \mathcal{A} \rightarrow \mathbb{R}^+$ represents the \textit{(state-action) transition probability density}. Let $\Gamma \in \mathcal{B}(\mathcal{S})$ be a Borel set of the state space, then the probability of transitioning into $\Gamma$ conditional on being in the state $s \in \mathcal{S}$ and having chosen action $a \in \mathcal{A}$ is given by 
\begin{equation*}
\mathbb{P}(s_{t+1} \in \Gamma \mid s_t=s, \, a_t=a ) = \int\limits_{\Gamma} p(s' | s, a) \mathrm{d} s' .
\end{equation*}
For more details about continuous MDP we recommend to look at \cite{VanHasselt2012}. For ease of notation, we will use the same symbol $p$ to denote the transition probability function and the transition probability density throughout the article.

\item and the \textit{reward function} $r: \mathcal{S} \times \mathcal{A} \rightarrow \mathbb{R}$ is the reward signal the agent will receive after being in state $s \in \mathcal{S}$ and having taken action $a \in \mathcal{A}$, i.e.
\begin{equation}
\label{eq: reward signal}
r_t \coloneqq r(s_t, a_t) .
\end{equation}
\end{itemize}
Formally, the tuple $(\mathcal{S}, \mathcal{A}, \mathbb{T}, p, r)$ defines a Markov decision process. A more detailed introduction to MDPs can be found in \cite{Puterman1994}.

\subsection{Reward and state-action transition probabilities}
\label{subsec: reward and transition probs}
If an action $a_t$ is chosen in a state $s_t$ at time $t$, two things happen. First, the agent receives a reward signal $r_{t}=r(s_t, a_t)$. Second, the agent transitions to the next state $s_{t+1}$ according to the transition probability function $s_{t+1} \sim p(\cdot |s_t, a_t)$. In the literature one can find reward functions that depend on the state, the action and the next state. In this case the order changes. First the agent moves to the next state, and second, the agent receives the reward signal. We will consider the first case throughout this paper. Note that both the (state-action) transition probability function and the reward function depend only on the current state of the agent and the action chosen in that state. This is sufficient to describe the dynamics of the agent, since we assume that the agent is Markov. Recall that the Markov decision framework can be generalised to non-Markovian dynamics; see, e.g., \cite{Puterman1994}. \\

The reward signal depends solely on the reward function. The reward function cannot be influenced by the agent and is considered to be part of the environment. We use the convention that the reward function is deterministic. However, for some environments the reward function may be described in a probabilistic way. In these cases it is useful to work with the so-called \textit{dynamics function} $p_\text{dyn}: \mathcal{S} \times \mathbb{R} \times \mathcal{S} \times \mathcal{A} \rightarrow [0, 1]$ which gives us the probability that the agent is in the next state and has received a certain reward conditional on a state-action pair \cite{Sutton2018}, i.e.
\begin{equation*}
p_\text{dyn}(s', r, s, a) = p_\text{dyn}(s', r | s, a) \coloneqq \mathbb{P}(s_{t+1}=s', r_{t}=r \mid s_t=s, a_t=a).
\end{equation*}
The \textit{expected reward function} $\mathcal{R}: \mathcal{S} \times \mathcal{A} \rightarrow \mathbb{R}$ provides us the expected reward signal conditioned on being in a state and taking an action
\begin{equation*}
\label{eq: exp reward function}
\mathcal{R}(s, a) \coloneqq \mathbb{E}[r_{t} \mid s_t=s, \, a_t=a ].
\end{equation*}

If the dynamics of the environment is deterministic, the state transition probability function is replaced by the so-called \textit{environment transition function}, denoted by $h:\mathcal{S} \times \mathcal{A}\rightarrow \mathcal{S}$. In this case the next state is given by $s_{t+1} = h(s_t, a_t)$. By introducing in the transition function a dependence on a random disturbance $\xi_t$ one can treat both stochastic and deterministic dynamical systems $s_{t+1} = h(s_t, a_t, \xi_t)$. To treat deterministic systems in this framework we just have to set $\xi$ to zero; see, e.g., \cite{Recht2019}.

\subsection{Policies}
\label{subsec: policies}
Policies are the most important part of reinforcement learning. A policy is a mapping that determines what action to take when the agent is in state $s_t$ at time $t$. This is why they are sometimes called the agent's brain. Policies can be either deterministic or probabilistic. 

A deterministic policy is a function from the state space into the action space. Deterministic policies are usually denoted by $\mu$, and the action for state $s_t$ is given by $a_t=\mu(s_t)$. A stochastic policy is a conditional probability distribution over the action space. Stochastic policies are usually denoted by $\pi$, and the new action for state $s_t$ can be computed by sampling from this conditional probability distribution. The action for $s_t$ is given by $a_t \sim \pi(\cdot|s_t)$. The deterministic policies can be seen as special cases of stochastic policies where the probability distributions over the action space are degenerate. \\

In the following sections we will show how the goal of reinforcement learning can be expressed in terms of finding policies that maximise the reward signal received at each time step.

\subsection{Trajectories, return and value functions}
\label{subsec: trajectories return and value functions}

Depending on the time horizon which is considered the reinforcement learning literature distinguishes between \textit{infinite horizon problems} and \textit{terminal problems}. We are only going to discuss terminal problems here since the SOC problem belongs to this class. For infinite time horizon problems a discounted factor has to be taken into account to make sure that the cumulative reward is finite. Details can be found in \cite{Sutton1999}. \\

By interacting with the environment the agent generates a trajectory $\tau$ which is a sequence of states ($s_.$), actions ($a_.$), and rewards ($r_.$)
\begin{equation}
\label{eq: trajectory}
\tau \coloneqq s_0,a_0,r_0,s_1,a_1,r_1, \ldots s_k, a_k, r_k, \dots.
\end{equation}
The initial state $s_0 \in \mathcal{S}$ is either sampled from a start state distribution $ s_0 \sim \rho_0$ or chosen to be constant $s_0 = s_{\text{init}}$, i.e $s_0 \sim \rho_0 = \delta_{s_{\text{init}}}$.  Recall that the state transitions only depend on the most recent state and action $s_{t+1} \sim p(\cdot| s_t, a_t)$.

The overall goal of the agent is to maximize the cumulative reward along a trajectory up to time $t$
\begin{equation}
\label{eq: cumulative reward}
G_t(\tau) \coloneqq r_{t} + \dots + r_T = \sum_{k=0}^{T-t} r_{t+k}
\end{equation}
by choosing the actions or policy optimally.

In order to estimate how well the agent performs starting in a given state and following a certain policy we need to define a performance measure. For each policy this performance measure is called \textit{value function} which is a function of the state only. It is defined as the expected return conditional on the agent starting at the state $s \in \mathcal{S}$
\begin{equation}
\label{eq: value function}
V^\pi(s) \coloneqq \mathbb{E} \big [ G_t(\tau) \mid s_t = s ; \pi \big],
\end{equation}
where the actions are chosen according to the policy $a_t \sim \pi(\cdot|s_t)$ for all $t \in \{0, \dots, T \}$. In optimal control theory the value function term is a synonym of the so-called \textit{optimal cost-to-go} which refers to the optimal value of the cost functional (with respect to all possible controls) conditional on the initial condition. In RL the term value function is used for all policies. When the policy is optimal we use the term \textit{optimal value function}. A crucial property of the value function is that it satisfies the following recursive relationship. For all $s \in \mathcal{S}$ and $t \in \{0, \dots, T -1 \}$ we have

\begin{align*}
V^\pi(s) 
&= \mathbb{E} \left[ \sum_{k=0}^{T-t} r_{t+k} \Bigm| s_t = s ; \pi \right] 
= \mathbb{E} \big [ r_{t} + \ldots + r_{T} \mid s_t = s ; \pi \big] \\
&= \mathbb{E} \left [ r_{t} + \mathbb{E} \left [\sum_{k=0}^{T-1-t} r_{t+1+k} \Bigm| s_{t+1} = s_{t+1} ; \pi \right] \Bigm| s_t= s ; \pi \right] =  \mathbb{E} \big [ r_{t} +  V^\pi(s_{t+1}) \mid s_t = s ; \pi \big ].
\end{align*}
In the reinforcement learning literature often stationary problems are considered; see, e.g., \cite{Sutton2018}. However, from a theoretical view point is possible to extend this to time dependent problems.
Moreover, the value function gives us a partial ordering over policies. We say that a policy $\pi$ is better than or equal to another policy $\pi'$, $\pi \geq \pi'$, if $V^\pi(s) \geq V^{\pi'}(s)$ for all states $s \in \mathcal{S}$. If a policy is better than or equal to all other policies, this is the optimal policy. The optimal policy is denoted as $\pi^*$ and the optimal value function is defined by
\begin{equation}
V^*(s) = \max\limits_\pi V^\pi(s), \quad \forall s.
\end{equation}

For a small, finite MDP, policy iteration (PI) strategies offer in general convergence to the optimal policy \cite{Sutton2018}. They interlude value iteration where the value function is estimated by using the Bellman equation iteratively and policy improvement. A direct way to find the optimal actions is given by solving the Bellman equation. For any MDP results about the existence of an optimal policy can be found in \cite{Puterman1994}.

Analogous to the motivation of the value function one might be interested in estimating how well the agent performs following a policy starting in a given state and taking a certain action. The performance measure for this is the so-called \textit{Q-value function}. The Q-value function for a certain policy is a function of state and action and is given by
\begin{equation}
\label{eq: q-value function}
Q^\pi(s, a) \coloneqq \mathbb{E} \left [ G_t(\tau) \mid s_t =s, a_t = a ; \pi \right] .
\end{equation}
It provides us with the quality of the chosen action with respect to the given state. Since the Q-value function is an extension of the value function it also satisfies the Bellman expectation equation
\begin{equation}
Q^\pi(s, a) =  \mathbb{E} \left [ r_{t} +  Q^\pi(s_{t+1}, a_{t+1}) \mid  s_t =s, a_t = a ; \pi \right ]
\end{equation}
and thus the recursive relationship.
Last, let us define the \textit{advantage function} with respect to a certain policy by
\begin{equation}
\label{eq: advantage function}
A^\pi(s, a) \coloneqq Q^\pi(s, a) - V^\pi(s),
\end{equation}
which tells us for each state-action pair the advantage of taking action $a$ in state $s$ with respect to the value function at $s$.

\subsection{RL as optimization problem}
\label{subsec: rl optimization problem}
Solving the Bellman optimality equations is often not an easy task because the transition function has to be known beforehand or learned. Even if this is the case it becomes an unfeasible assignment for high-dimensional problems due to the curse of dimensionality. In reinforcement learning one way to circumvent this is to formulate the RL problem as an optimization problem. Very generally this optimization problem can be stated as 
\begin{equation}
\label{eq: rl opt problem}
\pi^* = \argmax\limits_{\pi \in \Pi} J(\pi)
\end{equation}
where $\Pi$ is the space of policies we consider and $J$ is the RL objective function
\begin{subequations}
\begin{align}
J(\pi) &\coloneqq \mathbb{E}[G_0(\tau) \mid \pi] \label{eq: rl objective function} \\
\text{s.t.} \, s_0 \sim \rho_0, \, s_{t+1} &\sim p(\cdot| s_t, a_t), \, a_t \sim \pi(\cdot| s_t), \label{eq: rl dynamics} \,
\end{align}
\end{subequations}
for all $t \in \{0, \dots, T-1\}$ where $a_t$ is the action chosen at the state $s_t$ following policy $\pi$ and $s_{t+1}$ is the next state of the agent given by the dynamics of the system. The expectation in \eqref{eq: rl objective function} denotes the expected return along all trajectories following the policy $\pi$
\begin{equation*}
\mathbb{E}[G_0(\tau) | \pi] = \int_{\widetilde{\Omega}} \rho_\text{traj}(\tau | \pi) G_0(\tau) d\tau,
\end{equation*}
where $\widetilde{\Omega}$ denotes the space of all possible trajectories and $\rho_\text{traj}(\cdot | \pi)$ represents the probability distribution over the trajectories which follow the policy $\pi$
\begin{equation}
\label{eq: prob distr trajectories}
\rho_\text{traj}(\tau| \pi) = \rho_0(s_0) \prod_{t=0}^{T-1} p(s_{t+1}|s_t, a_t) \pi(a_t| s_t).
\end{equation}

In general it is numerically infeasible to optimize over a function space or a space of probability distributions. Hence, one often considers a parametric representation of a policy and tries to optimize with respect to the chosen parameters.

\subsection{Brief summary of RL algorithms}
\label{subsec: rl algorithms}

Over the last few years many different algorithms have been developed to solve reinforcement learning problems. For a good but non-exhaustive review we refer to \cite{Li2018} and the references therein. Many of the algorithms share a general framework, which is summarized in \Cref{alg: general rl}.

\begin{algorithm}
\caption{Main RL Online model-free algorithm}\label{alg: general rl}
\begin{algorithmic}[1]
\For {$\textit{trajectory}=1,2,\ldots$}
    \For {$\textit{time step}=1,2,\ldots,T$}
		\State Evaluate the dynamical system with the current policy $\pi$ and calculate the reward.
	\EndFor
\State Optimize the policy.
\EndFor
\end{algorithmic}
\end{algorithm}

The main difference between the methods is how the policy optimization is computed which depends on the underlying problem. The proposed methods can be distinguished if they are model based or model-free and if the policy optimization is done via Q-learning, policy gradient or a combination of both methods. \\

Let us first discuss the difference between model-based and model-free approaches. An algorithm is said to be model-based if the transition function is explicitly known (tractable to evaluate) point-wise. In this case we can just take the transition function and sample $s'$ from $p(\cdot|s,a)$. The known transition function contains a lot of information about the underlying dynamics and is therefore very useful for possible solution methods. A method is called \textit{model-free} when the transition function is not explicitly known. In this case the transition function cannot be used explicitly in the solution methods. The difference is related to the distinction between stochastic optimal control and reinforcement learning. In the case of stochastic optimal control the model is often known while reinforcement learning aims at a more general solution method that does not depend on the underlying dynamical system. \\ 

Let us briefly discuss the main underlying solution methods.

\subsubsection{Q-learning}
One of the first methods which was proposed was Q-learning \cite{Watkins1992}. The main idea is to approximate the state-action value function and use greedy policy iteration until convergence to the optimal policy \cite{Sutton2018}. The Bellman optimality equation
\begin{equation*}
\label{eq: q-value bellman opt equation}
Q^*(s, a) 
=  \mathbb{E} \Bigl[ r_t +  \max_{a' \in \mathcal{A}}Q^*(s_{t+1}, a') \bigm| s_t =s, a_t = a, \pi^* \Bigr]
\end{equation*}
provides a recursive formula for updating the Q-value function until convergence as long as some mild assumptions with respect to the learning rates and how often a state action pair is visited are satisfied (see \cite{Watkins1992, Sutton2018} for more details). If this is the case for any state $s \in \mathcal{S}$ the optimal policy is given by the action which maximizes the Q-value function
\begin{equation*}
\pi^*(s) = \argmax\limits_{a \in \mathcal{A}} Q^*(s, a).
\end{equation*}
However, the Q-learning algorithm can only be applied to problems with finite action spaces since taking the maximum over a continuous action space is not a feasible task. It works well as long as the state space is small and discrete. For large or continuous state spaces Deep Q-learning has been developed in \cite{Mnih2013} where the Q-value function is approximated by a parametric representation called the Deep Q-value Network (DQN). In addition Q-learning algorithms suffer from possible overestimation of the action values due to the fact that the maximum expected action value is approximated by the maximum action value which is biased. The idea of Double Q-learning attempts to overcome this problem by storing two Q-value functions. The method has been introduced for both tabular and function approximation settings \cite{VanHasselt2010, VanHasselt2015}.

An extension of the Q-learning idea for continuous action spaces is only possible by considering a separate policy parameterisation leading to an actor-critical setting which will be discussed in detail in \Cref{sec: model-free dpg}. We refer to \cite{Jia2022} and the references therein for a comprehensive study of Q-learning in continuous time.

\subsubsection{Policy gradient}
The main idea is to approximate the policy by a parametric representation and then solve the RL optimization problem \eqref{eq: rl opt problem} with respect to the parameters. In this approach one would like to use a direct optimization method like gradient descent. But one big challenge is the gradient calculation of the expectation, which is a research field by its own; see, e.g., \cite{Fournie1999}. One way to overcome this is to consider a parametrized stochastic policy. In this case the gradient calculation can be done explicitly and one can find a close formula. This was first presented in \cite{Williams1992} and the so-called \textit{policy gradient theorem} for stochastic policies was derived. The derived estimator is given by
\begin{equation}
\begin{split}
\label{eq: stochastic policy gradient}
\nabla_\theta J(\pi_\theta)
&= \mathbb{E} \Bigl[G_0(\tau) \nabla_\theta \log \rho_\text{traj}(\tau | \theta) \bigm| \pi_\theta \Bigr] \\
&= \mathbb{E} \Biggl[G_0(\tau) \Biggl(\sum\limits_{t=0}^{T-1} \nabla_\theta \log{\pi_\theta(s_t, a_t)}\Bigr) \Bigm| \pi_\theta \Biggl] \\
&= \mathbb{E}\left[\sum\limits_{t=0}^{T-1} \left(\sum\limits_{t'=0}^T r_{t'} \right)\nabla_\theta \log{\pi_\theta(s_t, a_t)} \Bigm| \pi_\theta \right]
\end{split}
\end{equation}
where one first uses the so-called log-derivative trick and then expands the expression for the gradient of the logarithm of the distribution of trajectories following the policy $\pi$. Notice that this gradient estimator does not depend on the transition probability density $p$ and hence it is considered to be model-free. \\

The methods presented here are the two main general ideas behind many variants of RL algorithms and were presented very early in the reinforcement learning community. Over the years, many drawbacks have been identified and extensions and combinations of solution methods have been proposed. Since policy evaluation is quite expensive, off-policy methods have been developed, i.e. trajectories simulated under a different policy are used for the current optimization step. This can be done, for example, by using an importance sampling approach; see, e.g., \cite{Schulman2017TRPO}. Methods that only use policies that have been sampled with the current policy are called online. A combination of Q-learning and policy gradient has been proposed to overcome the high variance of a pure policy gradient \cite{Schulman2018AC}. These ideas have been further developed and methods such as TRPO \cite{Schulman2017TRPO} and PPO \cite{Schulamn2017} have been proposed. \\

All of these algorithms have been developed for stochastic policies but methods for deterministic policies have also been derived. In the next section we will introduce the family of model-free deterministic policy gradient algorithms which provide a policy gradient without needing to know the model explicitly. 

The selection presented here is far from complete but a more detailed discussion of RL is beyond the scope of this article. For a more detailed overview we refer to \cite{Weng2018, Li2018, Sutton2018}. Most of the developed methods have a specific domain of application so one needs to carefully consider whether a method can be applied to a specific problem at hand.

\section{The SOC problem as RL formulation}
\label{sec: soc problem as rl}

In this section we show how the importance sampling SOC problem can be formulated as a reinforcement learning problem. First we show how to define an MDP for the stochastic control problem. This allows us to construct an RL environment based on the time-discretised stochastic optimal control problem. Then we compare the optimization approaches for both problems and argue that both formulations have a large overlap. After that we discuss how a previously presented method for the SOC problem can be categorised as a reinforcement learning algorithm. Finally we present a family of RL algorithms designed to deal with environments such as our stochastic optimal control problem.

\subsection{Importance sampling SOC problem as RL environment}
\label{subsec: is soc rl environment}
To show how the stochastic control problem can be interpreted as a reinforcement learning environment we consider for simplicity its corresponding time-discretized formulation. For the time-discretized dynamics of \eqref{eq: controlled sde} we use an Euler-Maruyama discretization of the SDE see e.g. \cite{Higham2001}
\begin{equation}
s_{t+1} = s_t + \left( -\nabla V(s_t)  + \sigma u(s_t)\right) \Delta t + \sigma \, \sqrt{\Delta t} \, \eta_{t+1}, \quad s_0 = x
\label{eq:soc_dis_dynamics} 
\end{equation}
where $s_t$ represents the time evolution of the controlled dynamics for a time step $\Delta t$. The term $\sqrt{\Delta t} \, \eta_{t+1}$ is the so-called Brownian increment where $\eta_{t+1} \sim \mathcal{N}(0,1)$ is a random number from a normal distribution. The time-discretized objective function is given by
\begin{equation}
J(u; x) \coloneqq \mathbb{E} \left[ g(s_{T_u}) + \sum \limits_{t=0}^{T_u-1} f(s_t) \Delta t + \frac{1}{2} \sum\limits_{t=0}^{T_u-1} |u(s_t)|^2 \Delta t \Bigm| s_0 = x \right]. \label{eq: soc dis objective function}
\end{equation}
The \textit{state space} consists of all possible states $s \in \mathbb{R}^d$ in the domain and is therefore infinite and continuous, $\mathcal{S} = \mathbb{R}^d$. The \textit{action space} is given by the space of the SOC controls i.e. $\mathcal{A} = \mathbb{R}^d$. Recall that the set of decision epochs $\mathbb{T}$ is a discrete set. As for the uncontrolled process we assume that the stopping time for the controlled process $T_u$ is a.s. finite and hence for any trajectory, the set of decision epochs $\mathbb{T}$ is finite. \\

To get a better understanding of the policy under consideration let us first have a look at the control used in \eqref{eq:soc_dis_dynamics}. The main idea to turn the SOC problem into an optimization problem is to use a Galerkin projection into a space of weighted approach functions. This was first proposed in the literature by \cite{Hartmann2012}. We will denote the approximation by $\texttt{u}_\theta$ and it is given by
\begin{equation}
\texttt{u}_\theta(x) = \sum_{m=0}^M \theta_m b_m (x)
\end{equation}
where $b(x)$ are some approximation functions such as radial basis functions, polynomials and $\theta $ is the weight vector. It is also possible to approximate the control with nonlinear function approximations such as neural networks; see, e.g., \cite{RiberaQuerRichter2022}. Regardless of the choice of approximation functions the resulting control is still deterministic which is usually called feedback form control in the optimal control literature. Due to the deterministic control the resulting policy is deterministic and the optimization will be over the weights of the ansatz function. \\

The time evolution for each control is given by the time-discretised dynamics \eqref{eq:soc_dis_dynamics} and this is sufficient to simulate trajectories. Moreover, for the considered SDE, the corresponding \textit{transition probability density} for the deterministic policy can even be explicitly derived from \eqref{eq:soc_dis_dynamics} (see \cite{RiberaQuerRichter2022} for details)
\begin{equation}
\label{eq: langevin transition density}
p(s_{t+1}|s_t, a_t) = \left(\frac{\beta}{4 \pi \Delta t}\right)^{d/2} \exp{\left(- \frac{\beta \Delta t}{4} \left| \frac{s_{t+1} - s_t}{\Delta t} + \nabla V(s_t)  - \sqrt{2 \beta^{-1}} a_t  \right|^2 \right)}.
\end{equation}

The \textit{reward function} is defined such that the corresponding return along a trajectory equals the negative term inside the expectation of the time-discretized cost functional \eqref{eq: soc dis objective function}. The reward signal at time step $t$ reads
\begin{equation}
\label{eq: langevin reward function}
r_t = r(s_t, a_t) \coloneqq \left\{
\begin{array}{ll}
- f(s_t) \Delta t - \frac{1}{2} |a_t|^2 \Delta t & \text{if} \quad s_t \notin \mathcal{T} \\
-g(s_t) & \text{if} \quad s_t \in \mathcal{T}.
\end{array}
\right.
\end{equation}
Notice that the reward signal is in general not sparse since the agent receives feedback at each time step but the choice of the running cost $f$ and the final cost $g$ can influence this statement. Moreover, the \textit{return} along a trajectory $\tau$ looks like
\begin{equation}
G_0(\tau) = - g(s_{T}) - \sum\limits_{t=0}^{T-1} f(s_t) \Delta t - \frac{1}{2} \sum\limits_{t=0}^{T-1} |a_t|^2 \Delta t.
\end{equation}

To do this we have defined an MDP for the importance sampling SOC problem. So we change the viewpoint of the SOC problem to the RL viewpoint. Next let us have a brief look at the two optimization problems for reinforcement learning and stochastic optimal control. 

\subsection{Comparison between the optimization approaches}
\label{subsec: soc vs rl}
Let us start by looking at the two optimization problems given in \eqref{eq: rl opt problem} and \eqref{eq: soc problem}. Both problems are stochastic optimization problems. The RL optimization problem maximises the expected return and the importance sampling SOC problem minimises the time-discretized cost functional conditional on the initial position. The RL optimization problem is a bit more general as the forward time evolution is stated in a very general setting. The forward trajectories are determined by the probability transition density and the stochastic policy while for the stochastic optimal control problem the time evolution is explicitly given by a controlled SDE \eqref{eq:soc_dis_dynamics} for a chosen control. Furthermore, in the RL framework, technical conditions are rarely imposed on the policy. In contrast in the SOC optimization formulation it is known that the optimal control is deterministic and must satisfy some technical assumptions as shown above. \\

The time-discretized optimization problem for the SOC problem is given by 
\begin{subequations} 
\begin{align}
J(\pi_\theta) &= \mathbb{E} \left[ g(s_{T}) + \sum \limits_{t=0}^{T-1} f(s_t) \Delta t + \frac{1}{2} \sum\limits_{t=0}^{T-1} |a_t|^2 \Delta t \Bigm| \pi_\theta \right] \rightarrow min\\
\text{s.t.} \quad s_{t+1} &= s_t + \left( -\nabla V(s_t)  + \sigma a_t \right) \Delta t + \sigma \, \sqrt{\Delta t} \, \eta_{t+1}, \quad s_0 = x, \quad a_t \sim \pi_\theta(\cdot|s_t) 
\end{align}
\end{subequations}
where $s_t$ is the time evolution of the controlled dynamics. Furthermore, we can turn the minimisation problem into a maximisation problem by multiplying the cost functional by $-1$ without making any further changes. Looking more closely at the reward signal for the SOC problem we see that the signal is a compromise between maximising the quantity of interest as much as possible which is given here by
\[
-g(s_T) - \sum\limits_{t=0}^{T -1} f(s_t) \Delta t,
\]

where the control term should not be too strong, i.e. 
\[
- \frac{1}{2} \sum\limits_{t=0}^{T - 1} |a_t|^2 \Delta t.
\]
So the last part of the reward can be seen as a kind of regularisation. 

In both cases the expectation is about trajectories. For the SOC case the probability density over the trajectories can be given in the same way as in \eqref{eq: prob distr trajectories}. Note that in this particular case the reward along a trajectory depends on the chosen policy. 
\subsection{Algorithms for SOC in a RL framework}
\label{subsec: algorithms for soc}

The first method proposed to solve the importance sampling SOC problem \cite{Hartmann2012} was a pure gradient descent but many other variants of this approach have been developed in the literature; see, e.g., \cite{Hartmann2017} and the references therein. Nevertheless most methods are based on the idea of gradient descent. We will show a derivation of the gradient descent method in the RL framework and show that this can be interpreted as a deterministic model-based policy gradient method. For the deterministic policy setting methods in reinforcement learning have also been proposed namely the model-free DPG family of algorithms and its subsequent variants. We will introduce the main idea of these methods and in the next section we will present an application of both algorithms.

\subsubsection{Model-based deterministic policy gradient}
Let us again understand the expectation of the objective function \eqref{eq: rl objective function} taken over the distribution of trajectories and proceed with the gradient computation as we did in \eqref{eq: stochastic policy gradient}.
First, notice that the probability distribution over the trajectories now reads
\begin{equation*}
\rho_\text{traj}(\tau| \mu_\theta) = \prod_{t=0}^{T-1} p(s_{t+1}|s_t, \mu_\theta(s_t))
\end{equation*}
and the derivative of the logarithm of the probability distribution with respect to the parameters changes to
\begin{align*}
\nabla_\theta \log\rho_\text{traj}(\tau | \mu_\theta)
&= \nabla_\theta \log{\left( \prod\limits_{t=0}^{T-1} p(s_{t+1} | s_t, \mu_\theta(s_t))\right)} \\
&= \sum\limits_{t=0}^{T-1} \nabla_\theta \log{(p(s_{t+1} | s_t, \mu_\theta(s_t)))} \\
&= \sum\limits_{t=0}^{T-1} \nabla_a \log p(s_{t+1} | s_t, a) \Big|_{a=\mu_\theta(s_t)} \cdot \nabla_\theta \mu_\theta(s_t).
\end{align*}
Moreover, the return along a trajectory now depends directly on the policy parameter $\theta$ 
\begin{equation*}
G_0(\tau; \theta) = \sum\limits_{t=0}^T r_t = \sum\limits_{t=0}^T r(s_t, \mu_\theta(s_t))
\end{equation*}
and therefore its derivative cannot be neglected
\begin{equation*}
\nabla_\theta G_0(\tau; \theta) = \sum\limits_{t=0}^T \nabla_\theta r(s_t, \mu_\theta(s_t)) = \sum\limits_{t=0}^T \nabla_a r(s_t, a) |_{a=\mu_\theta(s_t)} \nabla_\theta \mu_\theta(s_t).
\end{equation*}
By putting everything together the deterministic policy gradient reads
\begin{align}
\nabla_\theta J(\mu_\theta)
&= \mathbb{E} \Bigl[\nabla_\theta G_0(\tau; \theta) + G_0(\tau; \theta) \nabla_\theta \log\rho_\text{traj}(\tau | \mu_\theta) \bigm| \mu_\theta \Bigr] \\
&= \mathbb{E} \Biggl[\sum\limits_{t=0}^T \nabla_a r(s_t, a) |_{a=\mu_\theta(s_t)} \nabla_\theta \mu_\theta(s_t) \nonumber\\
& \quad\quad\quad + G_0(\tau; \theta) \left(\sum\limits_{t=0}^{T-1} \nabla_a \log p(s_{t+1} | s_t, a) \Big|_{a=\mu_\theta(s_t)} \cdot \nabla_\theta \mu_\theta(s_t) \right)  \Bigm| \mu_\theta \Biggr].
\end{align}
As mentioned by \cite{Peters2010} this gradient estimator requires the knowledge of the (state-action) transition density and the reward function. Hence, in contrast with its analog for stochastic policies it is a model-based method. Moreover, notice that this gradient estimator has the same form as \eqref{eq: stochastic policy gradient} with an extra term in the expectation corresponding to the gradient of the return. Thus an optimization method which uses this gradient estimator can be interpreted as a deterministic version of the well-known REINFORCE method. A different approach to derive the same gradient estimator for the continuous-time SOC problem can be found in \cite{Lie2021, RiberaQuerRichter2022}. \\

For the importance sampling RL environment we can provide a close expression for this gradient estimator since the model is known 
\begin{equation}
\label{eq:reinforce_grad}
\nabla_\theta J(\mu_\theta)
= \mathbb{E} \left[- \Delta t \, \sum\limits_{t=0}^{T -1} \mu_\theta(s_t) \cdot \nabla_\theta \mu_\theta(s_t) + G_0(\tau; \theta) \left(\sum\limits_{t=0}^{T -1} \eta_{t+1} \cdot \nabla_\theta \mu_\theta(s_t) \right) \Bigm| \mu_\theta\right]
\end{equation}

where the term $\nabla_a \log p(s_{t+1} | s_t, a)$ has been computed by using the Euler-Marujama time-discretized stochastic differential equation \eqref{eq:soc_dis_dynamics}. Such an approach has already been successfully implemented for the mentioned environment in the context of stochastic optimal control optimization. In the original work \cite{Hartmann2012} the deterministic policy i.e. the control is represented by a linear combination of Gaussian ansatz functions and later the approach is extended to a general deep representation of the control to avoid the curse of dimensionality \cite{RiberaQuerRichter2022}.

\subsubsection{Model-free deterministic policy gradient}
\label{sec: model-free dpg}
Due to the nature of the importance sampling SOC problem we are interested in exploring RL algorithms designed for continuous environments with optimal deterministic policies. The well-known deterministic policy gradient (DPG) algorithms presented in \cite{Silver2014} attempt to find a deterministic policy and are suitable for environments with continuous actions. The authors extend the \textit{policy gradient theorem} \cite{Sutton1999} to problem with deterministic policies and discounted infinite horizon. The gradient derived in this paper is given by
\begin{equation}
\label{eq: dpg}
\nabla_\theta J(\mu_\theta) = \mathbb{E}_{s \sim \rho^\mu} \bigl[\nabla_\theta \mu_\theta(s) \nabla_a Q^{\mu_\theta}(s, a)|_{a=\mu_\theta(s)} \bigr]
\end{equation}
where $\rho^\mu$ is meant to be the so-called ``(improper) discounted state distribution'' \cite{Sutton1999, Silver2014, Sutton2018}. The \textit{actor} representation $\mu_\theta$ approximates the deterministic policy and its parameters are updated by stochastic gradient ascent of \eqref{eq: dpg}. The \textit{critic} representation $Q_\omega$ approximates the Q-value function $Q^{\mu_\theta}$ and can be estimated by using policy evaluation algorithms. \\

Furthermore, the authors show that this deterministic policy gradient is the limiting case of its stochastic analogue as the policy variance tends to zero. This is shown by using a stochastic policy $\pi_{\mu,\sigma}$, where $\mu$ is a parameter for the mean and $\sigma$ controls the variance of the underlying probability distribution. To find a gradient estimator for deterministic policies, the log tick is applied to obtain an analytic expression for the gradient. The variance parameter is then set to $0$. Under certain assumptions, it can be shown that the following expression holds (see \cite{Silver2014})
\begin{equation}
\lim_{\sigma \rightarrow 0} \nabla_\theta J(\pi_{\mu_\theta,\sigma}) = \nabla_\theta J(\mu_\theta) .
\end{equation}
The underlying idea is that the probability distribution degenerates and the resulting degenerated probability distribution is only optimized in the mean parameter. 

Later this actor-critic approach was developed further resulting in a Deep Deterministic Policy Gradient (DDPG)\cite{Lillicrap2019} version of the algorithm where the policy and the Q-value function are approximated by deep neural networks.
In this work the main ideas of deep Q-learning \cite{Mnih2013} have been adapted to problems with continuous actions. The critic network is updated by minimising the following loss function known as the \textit{(mean squared) Bellman loss}
\begin{equation}
\label{eq: bellman loss}
\mathcal{L}_{BE}(Q_\omega^{\mu}) = \mathbb{E}_{(s, a, r, s', d) \sim \mathcal{D}} \Bigl[ \Bigl( Q_\omega(s, a) - \Bigl( r(s, a) + (1 - d) Q_\omega(s', \mu(s')\Bigr) \Bigr)^2 \Bigr]
\end{equation}
where $\mathcal{D}$ is a given data set of \textit{trajectory transitions} and $d$ tells us if the next state $s'$ is a terminal state i.e. $d = \mathbbm{1}_{\mathcal{T}}(s')$. The Bellmann loss tells us how close the approximation $Q_\omega$ is to satisfy the Bellmann equation in the given data set of transitions. \\

Further ideas behind the success of DQN algorithms \cite{Mnih2013} such as off-policy training with samples from a \textit{replay buffer} or the use of separate \textit{target networks} have been implemented to provide more stable and robust learning. The replay buffer consists of a finite set of trajectory transitions $(s_t, a_t, r_t, s_{t+1}, d)$ that are stored online or offline. When the replay buffer is full the oldest tuples are discarded. The transition records required in \eqref{eq: bellman loss} are randomly sampled from the replay buffer. The motivation for using target networks is to reduce the correlations between the action values $Q_\omega(s, a)$ and the corresponding targets $r + (1 -d) Q_\omega(s', \mu(s))$. A separate network is used for the targets, and its weights are updated slowly according to the original network. By forcing the target networks to update slowly, the stability of the algorithm is improved. \\

Finally, the work of \cite{Fujimoto2018} addresses how to deal with a possible overestimation of the value function for continuous action space problems with ideas from Double Q-learning. They introduce the idea of clipped actions as a regularisation technique for deep value learning. The idea behind this is that similar actions should have similar action values. The algorithm developed is the latest developed algorithm in the DPG family and is called Twin Delayed Deep Deterministic policy gradient (TD3).

In general, the DPG family of algorithms combines the two ideas of policy optimization and Q-learning methods. On the one hand it is a type of policy gradient algorithm because it uses a parametric representation of the policy and updates its parameters by a gradient ascent method. On the other hand the required gradient depends on the Q-value function that needs to be approximated.
Last let us conclude by highlighting that the resulting DPG algorithms are model-free, i.e. they do not require knowing the model transition density (see \cite{Silver2014, Lillicrap2019} for further details). 
    
\section{Examples: 1-dimensional double well potential}
\label{sec: numerical examples}

In this section we compare the use of two different algorithms to solve the RL optimization problem for the importance sampling environment introduced in \Cref{sec: soc problem as rl}. We consider the two main possible approaches for environments with continuous states and actions where the desired optimal policy is deterministic. We restrict ourselves to deterministic policy methods because we know from the PDE connection of the SOC problem that the optimal solution is deterministic. The first is an online model-based policy gradient method and the second is an offline model-free actor-critic method. \\

Throughout, we will consider the paradigmatic 1-dimensional double well potential ${V_\alpha: \mathcal{D}_\mathcal{S} \subset \mathcal{S} \rightarrow \mathbb{R}}$ given by
\begin{equation*}
{V_\alpha(s) = \alpha (s^2 - 1)^2}
\end{equation*}
where the parameter $\alpha$ is responsible for the height of the barrier and $\mathcal{D}_\mathcal{S}=[-2, 2]$ is the domain of the state space. We assume that the dynamics of the environment is governed by the transition probability density \eqref{eq: langevin transition density} with $\sigma(x)=\sqrt{2 \beta^{-1}}$ and time step $\Delta t=0.005$. The parameter $\alpha$ and the inverse temperature $\beta$ encode the metastability of the system. The height of the barrier is set to $\alpha=1.0$ so that the metastability is mainly influenced by the choice of $\beta$. To evaluate the performance of both algorithms we set a non-metastable setting $\beta=1$ and a more metastable one $\beta=4$. The initial value of the trajectories is set to $s_\text{init} = -1$ and the target set is chosen to be $\mathcal{T} = [1, \infty)$ for all experiments. We set $f=1$ and $g=0$ so that the quantity of interest \eqref{eq: quantity of interest} is reduced to $I(x_T) = \exp (-T)$.

\begin{figure}[htbp]
\centering
\begin{subfigure}{0.44\textwidth}
\includegraphics[width=\textwidth]{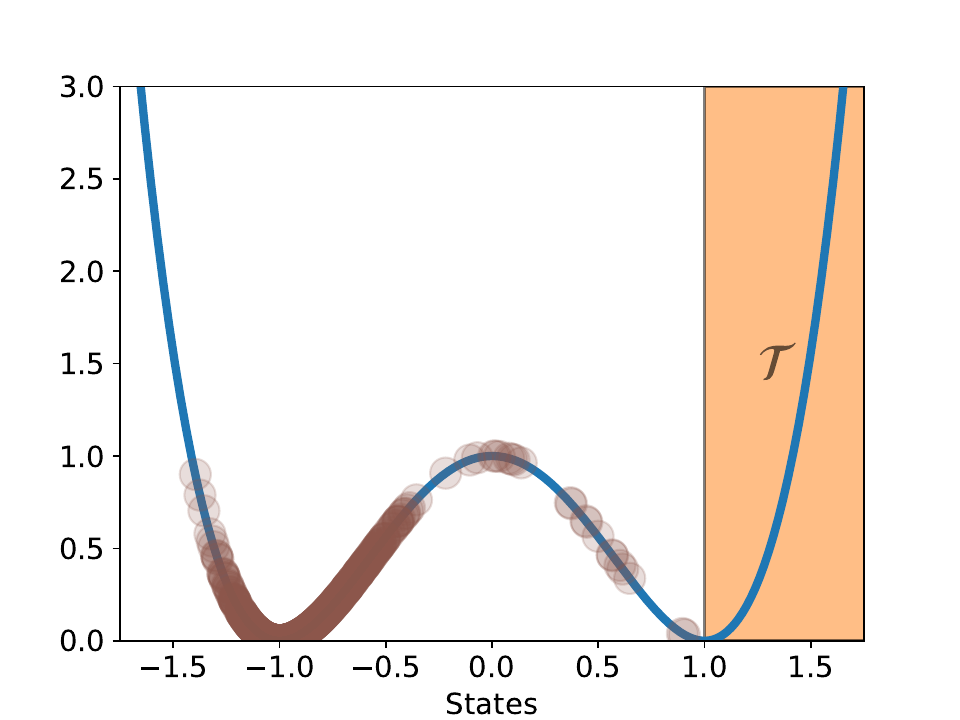}
\end{subfigure}
\begin{subfigure}{0.44\textwidth}
\includegraphics[width=\textwidth]{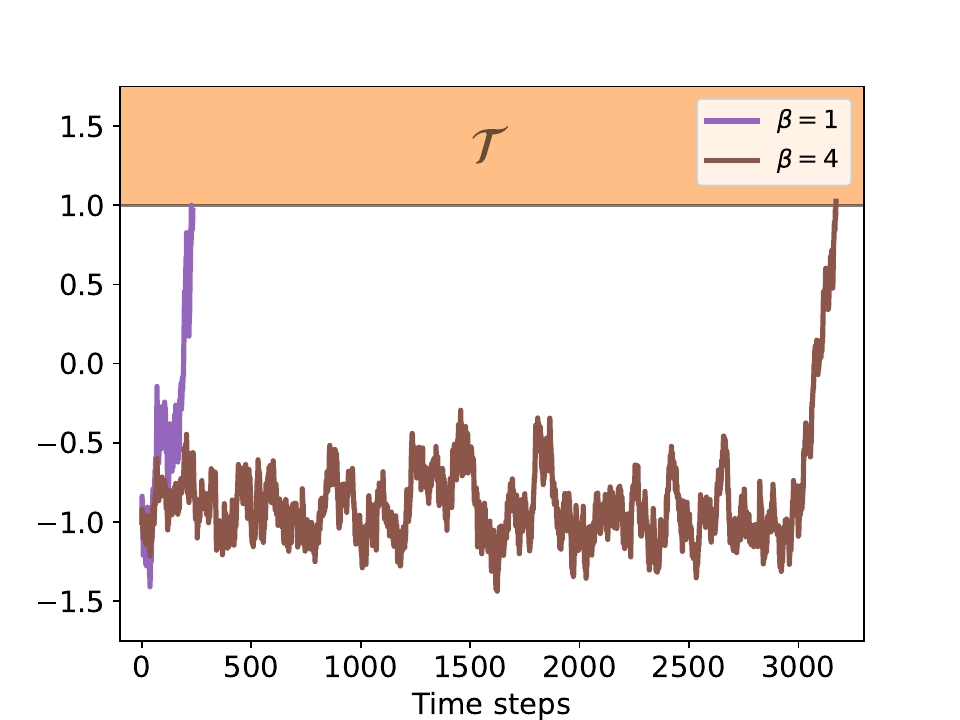}
\end{subfigure}
\caption{Trajectories following the not controlled policy for the two settings of study. The actions chosen along the trajectories are null. The trajectories are sampled starting at $s_\text{init} = -1$ until they arrive into the target set $\mathcal{T} = [1, \infty)$. Left panel: snapshots of the metastable trajectory $\beta=4$. Right panel: trajectory positions as a function of the time steps.}
\label{fig: setting}
\end{figure}

For the chosen environment the corresponding HJB equation \eqref{eq: hjb} is 1-dimensional and therefore we can compute reference solutions for optimal control by a finite difference method. We call this reference control the HJB policy. The construction of our environment is motivated by finding the optimal policy with the intention of minimising the variance of the original importance sampling estimator \eqref{eq: is estimator}. However, we do not focus on the performance of such variance reduction as it is sufficient to measure how close the approximated policy is to the optimal solution. To do this we track a $L^2$-type error of our deterministic policy approximations along the trajectories. Specifically, we define the empirical $L^2$ error along an ensemble of sampled trajectories with an arbitrary deterministic policy $\mu$ by

\begin{equation}
L^2(\mu) \coloneqq \frac{1}{K_\text{test}}\sum\limits_{k=1}^{K_\text{test}} \Biggl( \sum\limits_{t=0}^{T^{(k)}} |\mu - \mu_\text{HJB}|^2 (s_t^{(k)}) \Delta t \Biggr),
\end{equation}
where $K_\text{test}$ is the test batch size and the superscript $(k)$ denotes the index of each trajectory. This quantity tells us how close the policy $\mu$ is to the reference policy $\mu_\text{HJB}$ along the trajectories of the ensemble. For all experiments, the test batch size is chosen to be ${K_\text{test}=10^3}$. \newpage

To approximate the optimal policy and the Q-value function we use \textit{feed-forward neural networks}. These are essentially compositions of affine-linear maps and nonlinear activation functions which show remarkable approximation properties even in high dimensions. Let $d_\text{in}, d_\text{out} \in \mathbb{N}^+$ be the input and output dimensions of the feed-forward network $\varphi_\theta: \mathbb{R}^{d_\text{in}} \rightarrow \mathbb{R}^{d_\text{out}}$ which is defined recursively as follows
\begin{equation}
\label{eq: feed-foreward nn}
\varphi_\theta(x) = A_L \rho(A_{L-1} \rho(\cdots  \rho(A_1 x + b_ 1) \cdots) + b_{L-1}) + b_L,
\end{equation}
where $L$ is the number of layers $d_0 = d_\text{in}, d_L = d_\text{out}$, $A_l \in \mathbb{R}^{d_{l} \times d_{l-1}}$ and $b_l \in \mathbb{R}^{d_l}, 1 \le l \le L$ are the weights and the bias for each layer and $\rho: \mathbb{R} \to \mathbb{R}$ is a nonlinear activation function applied componentwise. The collection of matrices $A_l$ and vectors $b_l$ contains the learnable parameters $\theta \in \mathbb{R}^p$. For all experiments we consider the policy representation $\mu_\theta = \varphi_\theta$ where $d_\text{in} = d_\text{out} = 1$ and the Q-value representation $Q_\omega = \varphi_\omega$ where $d_\text{in} = 2$ and $d_\text{out} = 1$ with both $L=3$ layers, dimension of the hidden layers $d_1 = d_2 = 32$, and the activation function $\rho(x) = \tanh(x)$. To ensure that the initial output of the networks is close to zero the final layer weights and biases are initialised by sampling from the following uniform distributions $\mathcal{U}(-10^{-2}, 10^{-2})$, $\mathcal{U}(-10^{-3}, 10^{-3})$ respectively.

We repeat all our experiments several times with different random seeds to ensure generalisability. Each experiment requires only one CPU core, and the maximum value of allocated memory, is set to 1 GB unless otherwise stated.

\subsection{Model-based deterministic REINFORCE}
First we present the results of the deterministic model-based version of the REINFORCE algorithm for the 1-dimensional environment described above. We consider an online based implementation where the batch of trajectories is not reused after each gradient step. We summarise this method in \Cref{alg: det reinforce}.

\begin{algorithm}[H]
\caption{Model-based deterministic policy REINFORCE} \label{alg: det reinforce}
\begin{algorithmic} [1]
\State Initialize deterministic policy $\mu_\theta$.
\State Choose a batch size $K$, a gradient based optimization algorithm, a corresponding learning rate $\lambda > 0$, a time step size $\Delta t$ and a stopping criterion.
\Repeat
\State Simulate $K$ trajectories by running the policy in the environment's dynamics.
\State Estimate the policy gradient $\nabla_\theta J(\mu_\theta)$ by
\begin{equation*}
\frac{1}{K} \sum\limits_{k=1}^K \sum\limits_{t=0}^{T^{(k)} -1} \Bigl(- \Delta t \, \mu_\theta(s_t^{(k)}) \cdot \nabla_\theta \mu_\theta(s_t^{(k)}) + G_0(\tau; \theta) \eta_{t+1} \cdot \nabla_\theta \mu_\theta(s_t^{(k)}) \Bigr).
\end{equation*}
\State Update the parameters $\theta$ based on the optimization algorithm.
\Until{stopping criterion is fulfilled.}
\end{algorithmic}
\label{alg: reinforce deterministic policy}
\end{algorithm}

For both metastable and non-metastable settings we consider two different batch sizes $K = \{1, 10^3\}$ and use the Adam gradient based optimization algorithm \cite{Kingma2014} with corresponding learning rates $\lambda = \{5 \cdot 10^{-5}, 5 \cdot 10^{-4}\}$ respectively.
The stopping criteria is set to be a fixed number of gradient steps $N = 10^4$ and the approximated policy is tested every $100$ gradient updates. \\

\begin{figure}[htbp]
\centering
\begin{subfigure}{0.44\textwidth}
\includegraphics[width=\textwidth]{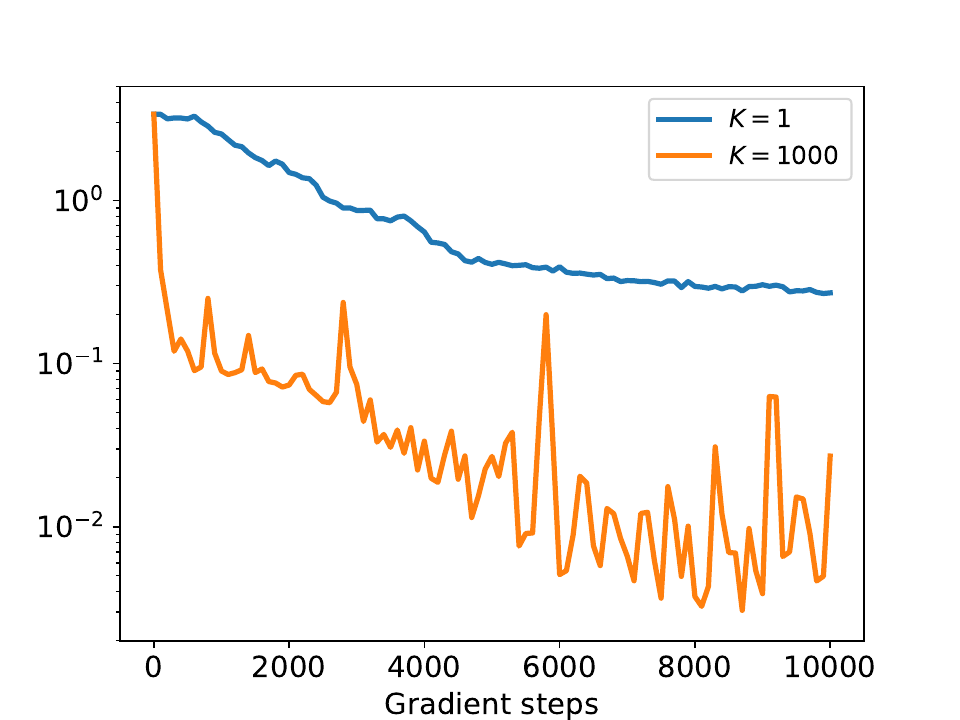}
\end{subfigure}
\begin{subfigure}{0.44\textwidth}
\includegraphics[width=\textwidth]{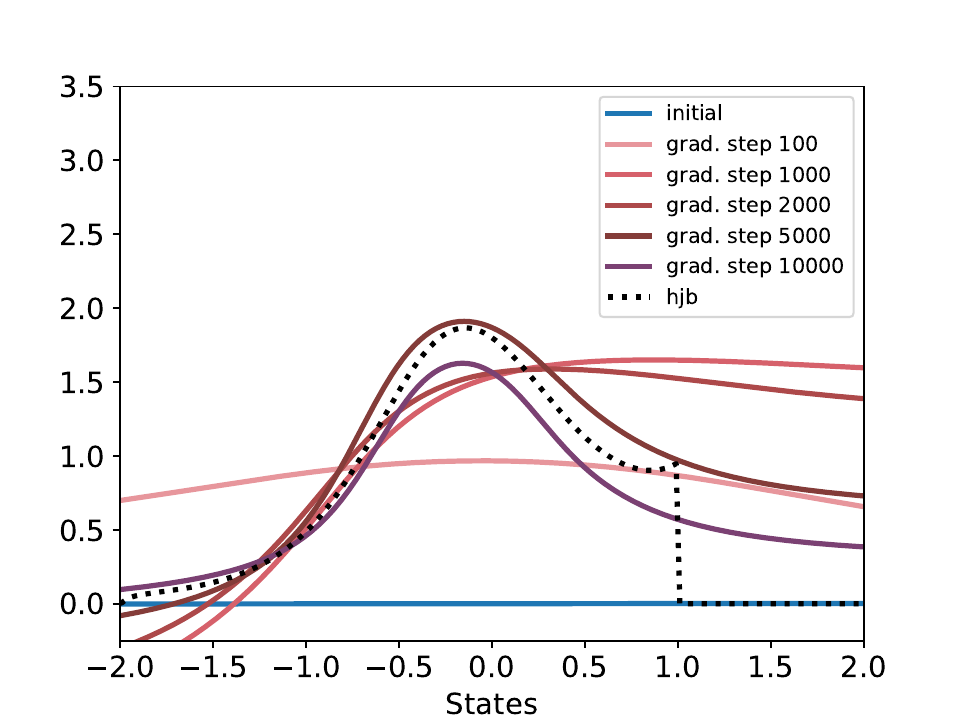}
\end{subfigure}
\caption{Left panel: estimation of $L^2(\mu_\theta)$ at each gradient step for the non-metastable setting $\beta=1$. Right panel: approximated policy for different gradient updates for the batch of trajectories case ($K=10^3$).}
\label{fig: det reinforce not metastable}
\end{figure}

\begin{figure}[htbp]
\centering
\begin{subfigure}{0.44\textwidth}
\includegraphics[width=\textwidth]{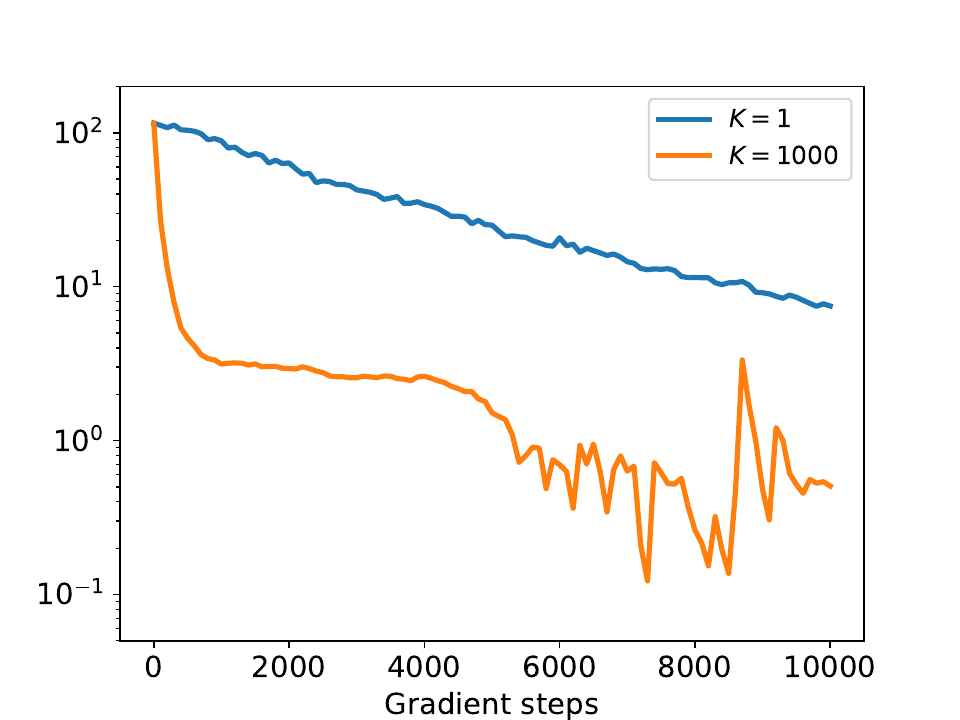}
\end{subfigure}
\begin{subfigure}{0.44\textwidth}
\includegraphics[width=\textwidth]{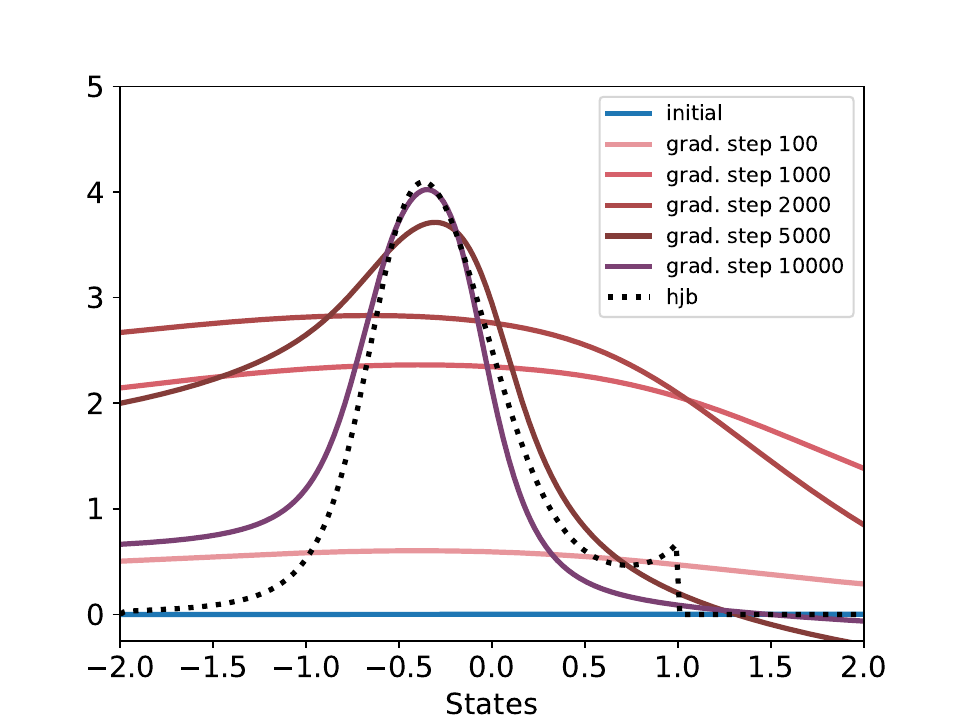}
\end{subfigure}
\caption{Left panel: estimation of $L^2(\mu_\theta)$ as a function of the gradient steps for the metastable setting $\beta=4$. Right panel: approximated policy for different gradient updates.}
\label{fig: det reinforce metastable}
\end{figure}

\Cref{fig: det reinforce not metastable} and \Cref{fig: det reinforce metastable} show the $L^2$ empirical error as a function of the gradient updates and the policy approximation at different gradient steps for the two different problem settings. In the non-metastable case we can see for the experiment with batch size $K=10^3$ that the policy approximation agrees well with the reference control already after $\sim 5 \, 000$ gradient steps. On the other hand with only one trajectory one has to rely on a lower learning rate and therefore learning is much slower. In the more metastable scenario we can see that the policy approximation for $K=10^3$ is not as close to the HJB policy as in the less metastable case. One can see that after the same number of gradient steps the $L^2$ error differs by more than an order of magnitude. Moreover, for the non-metastable setting with $K=10^3$ we need to increase the maximum value of allocated memory to 8GB. For the metastable setting we end up allocating 4GB for the one trajectory case $K=1$ and around 100GB for the batch case $K=10^3$. \newpage

\subsection{Model-free deterministic policy gradient}
Next we present the application of a model-free DPG method. In particular we will implement the TD3 variant of the DDPG method introduced at the end of \Cref{sec: model-free dpg}. First we consider a replay buffer and use separate target networks with slow updating so that the learning of the Q-value function is stabilised. We consider two critic networks to avoid overestimation of the value function. The replay buffer helps to reduce the amount of data that needs to be generated for each gradient estimation step. It uses caching of trajectories and random sampling of these cached trajectories for gradient estimation. The method is therefore considered an offline algorithm since it does not use the current control to estimate the gradient of the actor and critic networks. We summarise this method in \Cref{alg: td3}. \\

\begin{figure}[ht]
\centering
 \begin{subfigure}{0.44\textwidth}
\includegraphics[width=\textwidth]{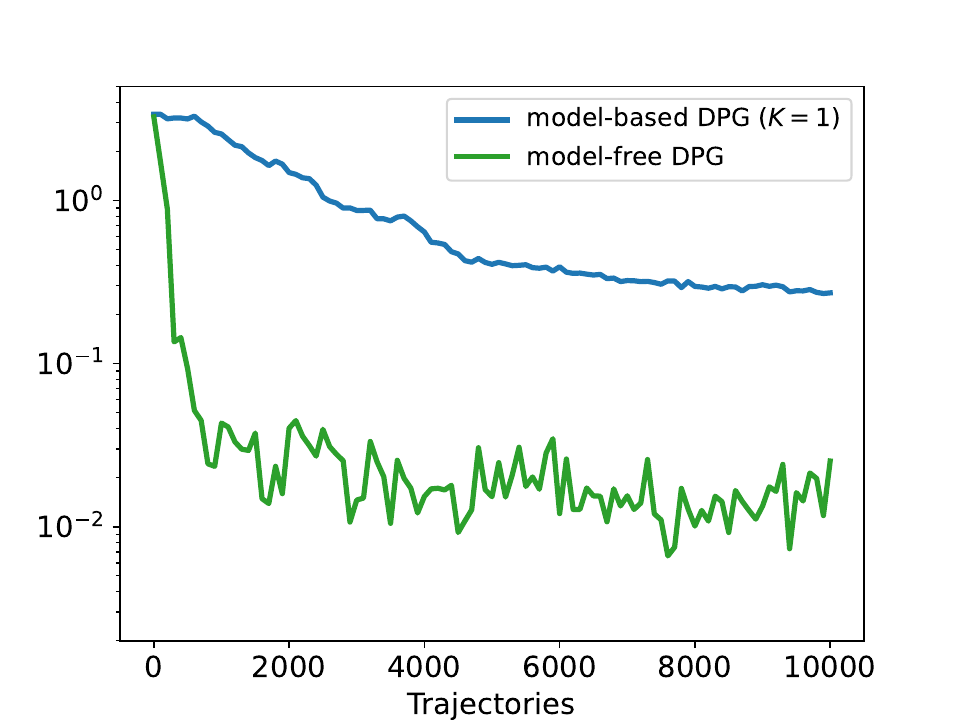}
\end{subfigure}
\begin{subfigure}{0.44\textwidth}
\includegraphics[width=\textwidth]{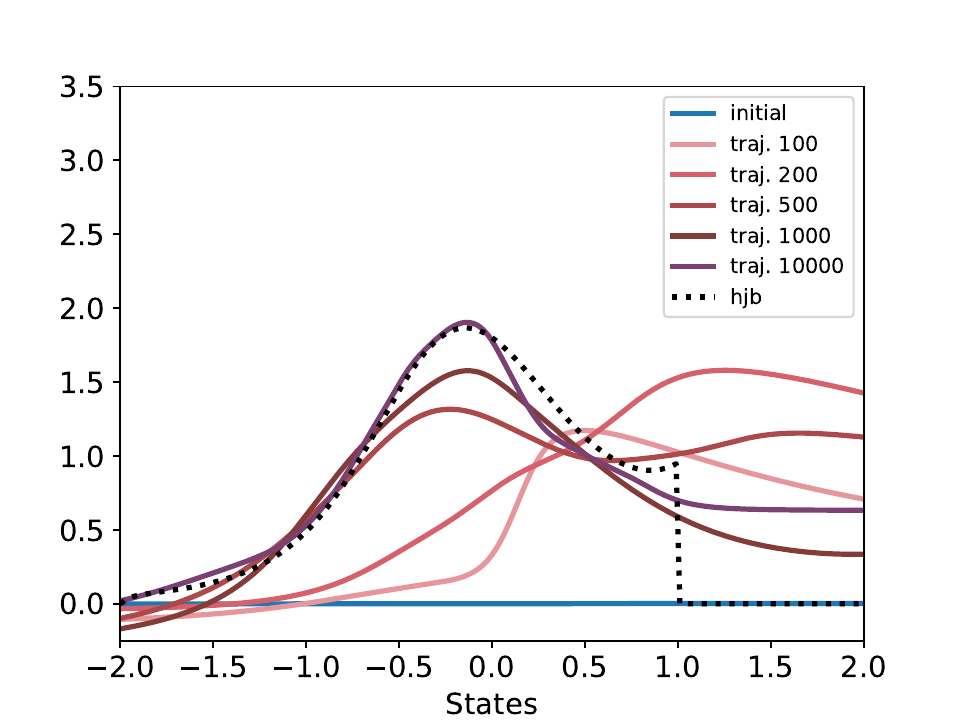}
\end{subfigure}
\caption{Left panel: estimation of $L^2(\mu_\theta)$ as a function of the trajectories for the non-metastable setting $\beta=1$. Right panel: approximated policy by the actor model after different trajectories.}
\label{fig: td3 not metastable}
\end{figure}

\begin{figure}[ht]
\centering
\begin{subfigure}{0.44\textwidth}
\includegraphics[width=\textwidth]{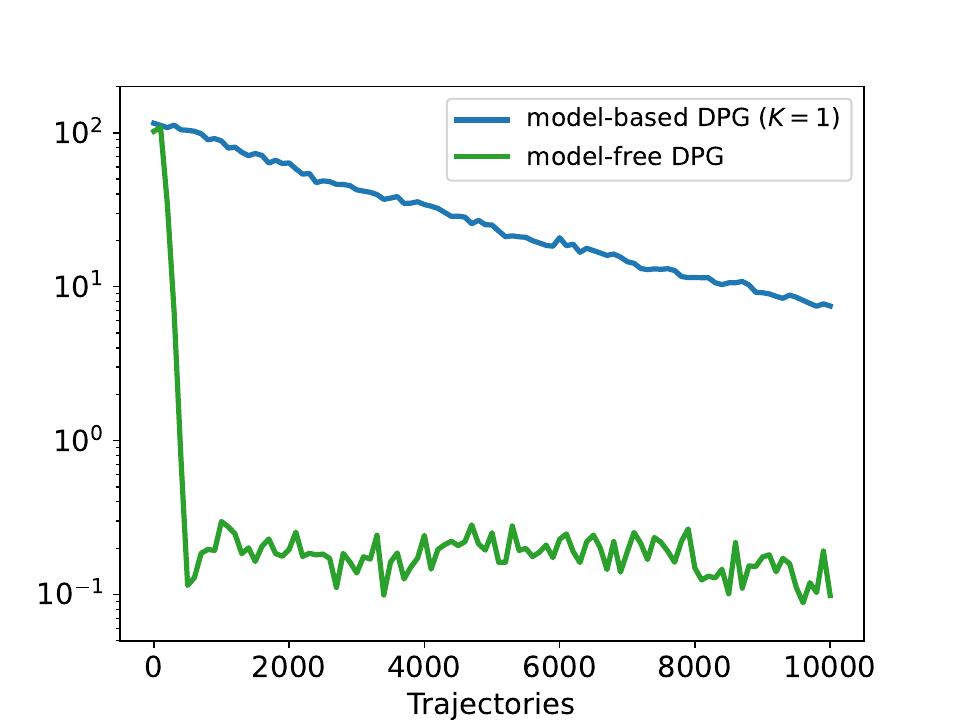}
\end{subfigure}
\begin{subfigure}{0.44\textwidth}
\includegraphics[width=\textwidth]{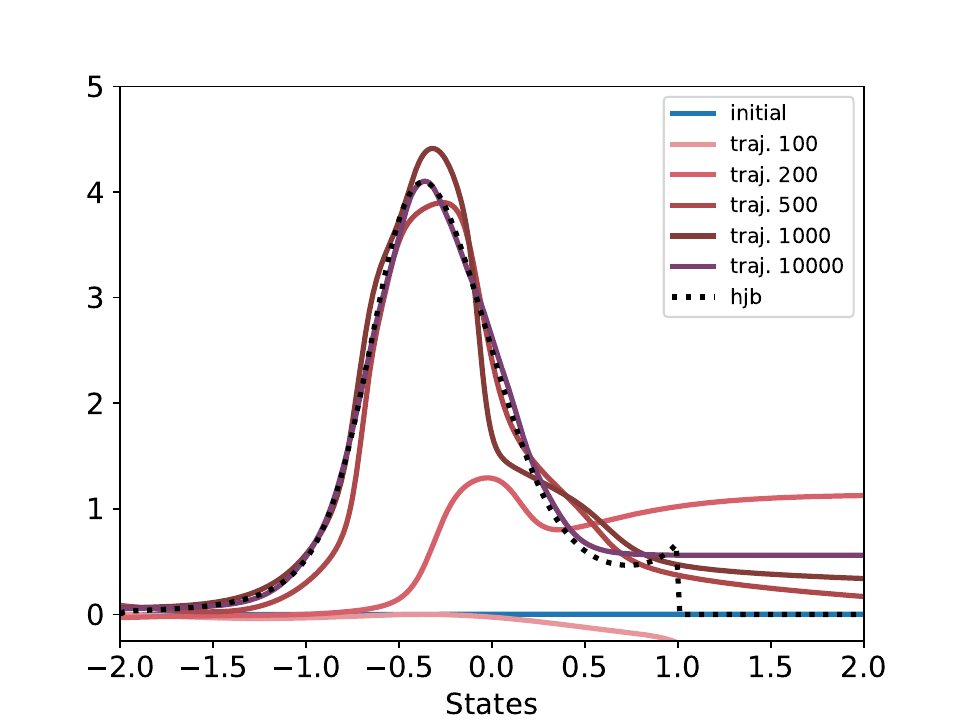}
\end{subfigure}
\caption{Left panel: estimation of $L^2(\mu_\theta)$ as a function of the trajectories for the metastable setting $\beta=4$. Right panel: approximated policy by the actor model after different trajectories.}
\label{fig: td3 metastable}
\end{figure}

\Cref{fig: td3 not metastable} and \Cref{fig: td3 metastable} show the evolution of the $L^2$ error as a function of the sampled trajectories and the policy approximation at the end of different trajectories for the two chosen settings. The $L^2$ error is compared with the model-based approach with one trajectory ($K=1$). For the non-metastable settings we can see that the policy approximation agrees well with the reference policy. We observe that the learning is quite fast. After around $2 \, 000$ the $L^2$ error is already smaller than $10^{-1}$, but after this point the error does not decrease any further. Moreover, we observe that the model-free method learns faster than the model-based method in terms of generated data  i.e. trajectories sampled. For the more metastable scenario we can see a similar pattern. Despite the metastability of the system the method learns quite fast and even manages to achieve a lower $L^2$ error than the model-based approach. However, it seems that the metastability does affect the stability of the method. Unfortunately this unstable behaviour is observed for other choices of the hyperparameters of the algorithm, where the $L^2$ error can even blow up.

\begin{figure}[ht]
\centering
\begin{subfigure}{0.44\textwidth}
\includegraphics[width=\textwidth]{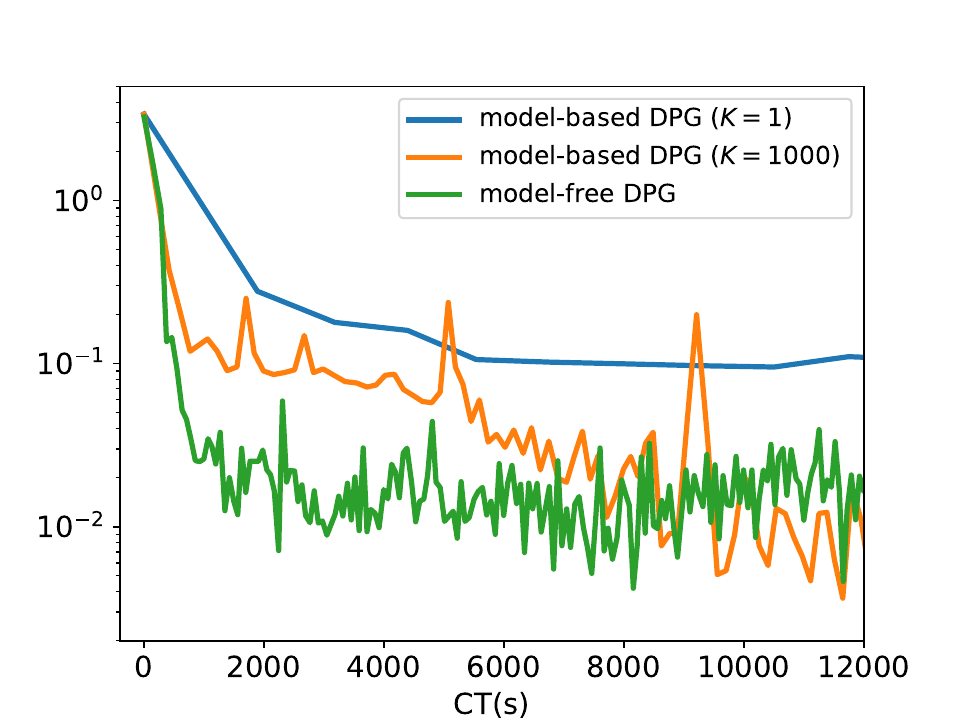}
\end{subfigure}
\begin{subfigure}{0.44\textwidth}
\includegraphics[width=\textwidth]{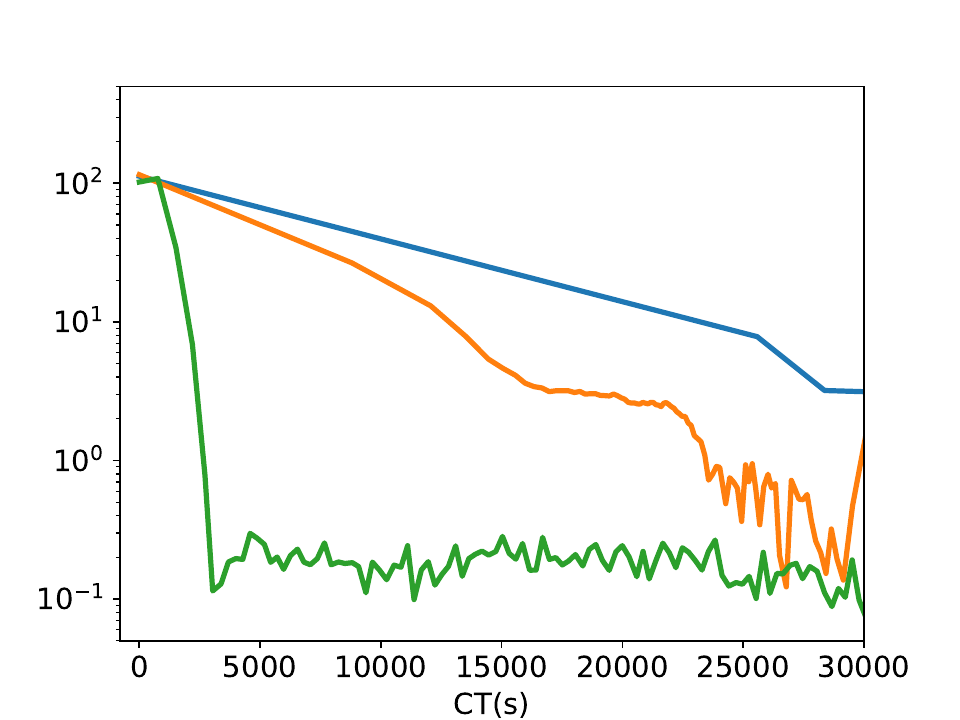}
\end{subfigure}
\caption{Estimation of $L^2(\mu_\theta)$ as a function of the computation time for the two considered methods. Left panel: non-metastable setting $\beta=1$. Right panel: metastable setting $\beta=4$.}
\label{fig: cts l2 errors}
\end{figure}

Finally in \Cref{fig: cts l2 errors} we compare the $L^2$ error as a function of the computation time for the two considered methods. We see that the TD3 method learns a decent control much faster than the model-based approach, especially in the metastable setting. However, we observe that at a certain point the $L^2$-error stops decreasing. In contrast to that, for the model-based method learning is much slower but there is a steady decrease in the $L^2$-error.

\begin{algorithm}[H]
\caption{Model-free deterministic Policy Gradient (TD3)} \label{alg: td3}
\begin{algorithmic} [1]
\State Initialize actor network $\mu_\theta$ and critic networks $Q_{\omega_1}$ and $Q_{\omega_2}$.
\State Initialize corresponding target networks: $\theta' \leftarrow \theta$, ${\omega'}_1 \leftarrow \omega_1$, ${\omega'}_2 \leftarrow \omega_2$ and choose $\rho_\text{p} \in (0, 1)$.
\State Initialize replay buffer $\mathcal{R}$.
\State Choose a batch size $K$, a gradient based optimization algorithm and a corresponding learning rate $\lambda_\text{actor}, \lambda_\text{critic} > 0$ for both optimization procedures, a time step size $\Delta t$ and a stopping criterion.
\State Choose standard deviation exploration noise $\sigma_\text{expl}$ and lower and upper action bounds $a_\text{low}, a_\text{high}$.
\Repeat
\State Select clipped action and step the environment dynamics forward.
\begin{equation*}
a= \text{clip}(\mu_\theta(s) + \epsilon, a_\text{low}, a_\text{high}), \quad \epsilon \sim \mathcal{N}(0, \sigma_\text{expl}).
\end{equation*}
\State Observe next state $s'$, reward $r$, and done signal $d$ and store the tuple $(s, a, r, s', d)$ in the replay buffer.
\If{$s'$ is terminal} 
\State Reset trajectory. 
\EndIf
\For{$j$ in range(\textit{update frequency})}
\State Sample batch $\mathcal{B} = \{(s^{(k)}, a^{(k)}, r^{(k)}, {s'}^{(k)}, d^{(k)})\}_{k=1}^K$ from replay buffer.
\State Compute targets (Clipped Double Q-learning and policy smoothing).
\begin{equation*}
y(r, s', d) = r + (1-d) \min\limits_{i=1, 2} \left\{ Q_{{\omega'}_i}(s', \tilde{a}) \right\}, \quad \tilde{a}= \text{clip}(\mu_{\theta'}(s) + \epsilon, a_\text{low}, a_\text{high}), \quad \epsilon \sim \mathcal{N}(0, \sigma_\text{target}). 
\end{equation*}

\State Estimate critic gradient $\nabla_\omega L(Q_\omega^{\mu_\theta})$ by
\begin{equation*}
\nabla_{\omega_i} \Biggl(\frac{1}{K} \sum\limits_{k=1}^K \left(Q_{\omega_i}(s^{(k)}, a^{(k)}) - y(r^{(k)}, {s'}^{(k)}, d^{(k)}) \right)^2 \Biggr), \quad \text{for} \,\, i=1,2.
\end{equation*}
\State Update the critic parameters $\omega_i$ based on the optimization algorithm.
\If{$j \text{ mod } \textit{policy delay frequency} = 0$}
\State Estimate actor gradient $\nabla_\theta J(\mu_\theta)$ by
\begin{equation*}
\nabla_{\theta} \Biggl(\frac{1}{K} \sum\limits_{k=1}^K Q_{\omega_1}(s^{(k)}, \mu_\theta(s^{(k)})) \Biggr).
\end{equation*}
\State Update the actor parameters $\theta$ based on the optimization algorithm.
\State Update target networks softly: 
\begin{equation*}
\theta' \leftarrow \rho_\text{p} \theta' + (1 - \rho_\text{p})\theta, \quad 
{\omega'}_i \leftarrow \rho_\text{p} {\omega'}_i + (1 - \rho_\text{p})\omega_i, \quad \text{for} \, \, i=1, 2.
\end{equation*}
\EndIf
\EndFor
\Until{stopping criterion is fulfilled.}
\end{algorithmic}
\end{algorithm} \newpage

\subsection{Discussion}
Let us now compare the two different methods used above. First we focus on the ingredients needed for each application. The model-based deterministic approach requires the model to be known. Without knowing the transition probability density this approach is not possible. For our importance sampling application with damped Langevin dynamics the transition probability density can be approximated after time discretization. However, we may be interested in general diffusion processes where this information is not given or cannot be trusted. Alternatively the model could be learned which is indeed a current area of research in model-based RL. On the other hand, the model-free alternative only requires knowledge of the reward function which is always the case for the importance sampling problem. \\

\begin{figure}[ht]
\centering
\begin{subfigure}{0.44\textwidth}
\includegraphics[width=\textwidth]{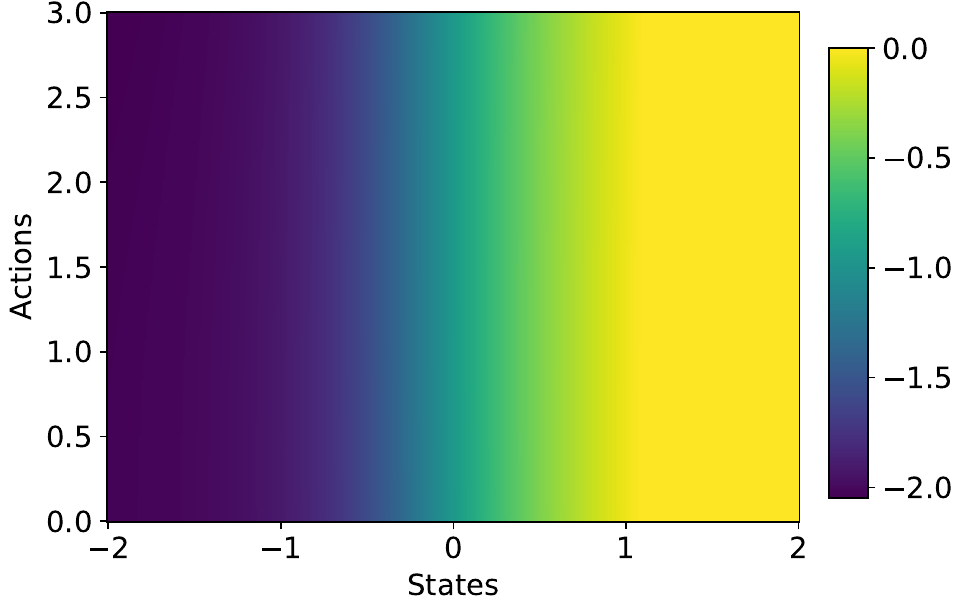}
\end{subfigure}
\begin{subfigure}{0.44\textwidth}
\includegraphics[width=\textwidth]{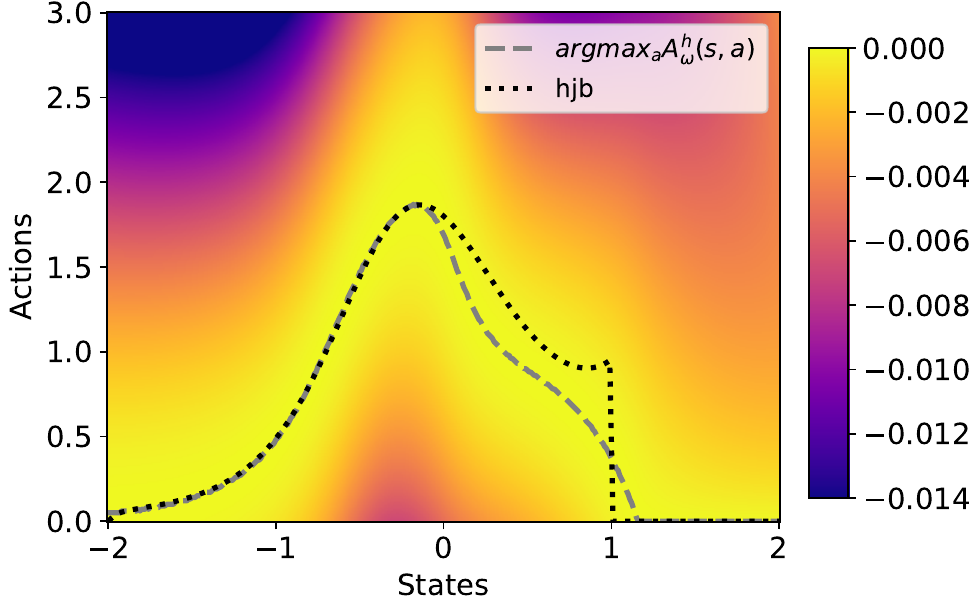}
\end{subfigure}
\begin{subfigure}{0.44\textwidth}
\includegraphics[width=\textwidth]{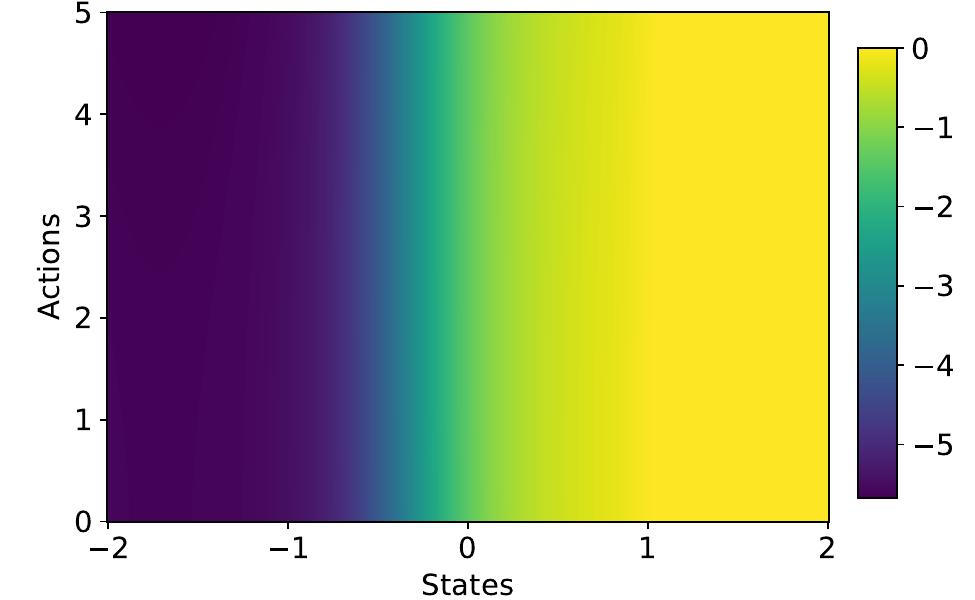}
\end{subfigure}
\begin{subfigure}{0.44\textwidth}
\includegraphics[width=\textwidth]{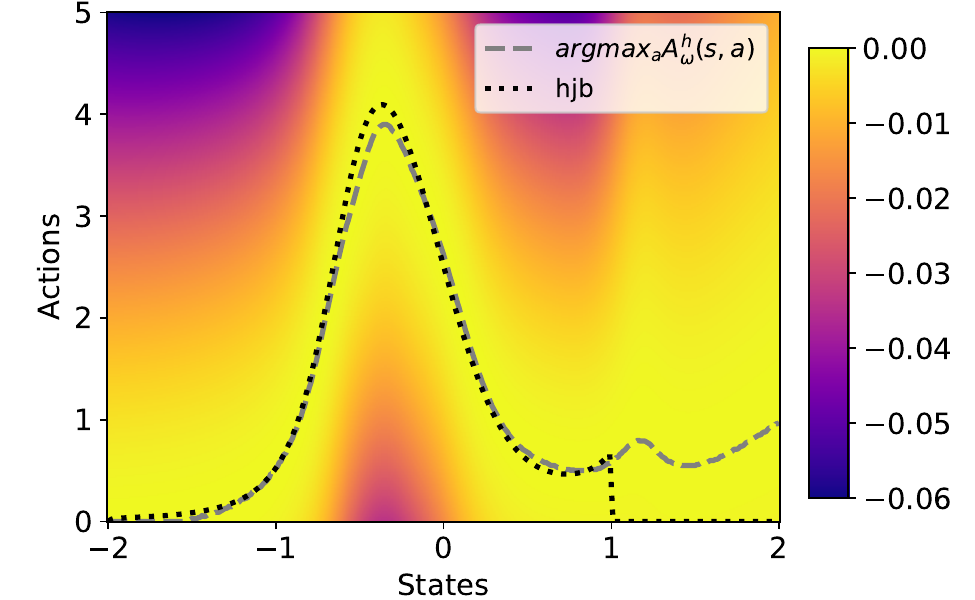}
\end{subfigure}
\caption{Critic models after the last trajectory for both settings. Left panels: approximated Q-value function $Q_\omega$. Right panels: advantage function after action space discretization $A_\omega^h(s, a) = Q_\omega(s, a) - \max\limits_{a \in \mathcal{A}_h} Q_\omega(s, a)$ and resulting greedy policy (grey dashed).}
\label{fig: td3 q-value and advantage function}
\end{figure}

However, this model-free approach has its implications. The method relies on a good approximation of the Q-value function especially along the action axis because this determines the direction of the gradient that the actor will follow. \Cref{fig: td3 q-value and advantage function} shows the approximated Q-value and the advantage function in the space-action discretized grid for both settings. Note that the shape of the Q-value function has the following property: the difference between the Q-values for a given action along the state axis is orders of magnitude larger than the difference between the Q-values for a given state along the action axis. This difference may be the cause of the observed instabilities in the model-free TD3 method. This problem is not specific to our importance sampling application and has been addressed in the field of Advantage Learning; see, e.g., \cite{Baird1998}. Advantage learning is an alternative approach to Q-learning where the advantage function is learned instead of the Q-value function. For future work it may therefore be interesting to exploit Dueling Network Architecture approaches \cite{Wang2015} where two separate estimators are maintained: one for the value function and one for the advantage function. \\

Regarding the performance of both approaches we have seen that the model-based method gets nearer to the optimal solution in the non-metastable settings. For high metastable dynamics this approach suffers from long running times and a high variance of the gradient estimator \cite{RiberaQuerRichter2022}. 
In our experimental analysis we observed a significant advantage of the TD3 algorithm over the deterministic reinforce algorithm in terms of learning a reasonable control faster. Our experiments suggest several reasons for this superiority. The TD3 algorithm performs a notably higher number of gradient steps per episode compared to the deterministic reinforce algorithm which relies on complete trajectory sampling before each gradient update, making bootstrapping impractical. This difference in gradient steps necessitates the deterministic gradient method to allocate considerably more memory for each update due to the extended length of trajectories. To ensure a fair comparison we conducted experiments for the case K = 1, allowing more gradient steps per data generated. The faster convergence observed in this scenario suggests that TD3 particularly benefits from the increased number of gradient steps especially in handling metastable problems.
Another critical aspect contributing to TD3's effectiveness is its reliance on accurate Q-value function approximation. When the Q-value function is well-approximated TD3's gradient updates effectively guide the algorithm towards the optimal policy without requiring complete trajectory sampling, unlike the deterministic gradient estimator which lacks this correction term. Furthermore, this advantage enables TD3 to pursue off-policy learning enhancing its overall efficiency and adaptability.
Moreover, TD3's third advantage lies in its integration of exploration mechanisms a feature absent in traditional gradient-based methods. By actively exploring the environment, TD3 efficiently uncovers novel and potentially rewarding state-action trajectories resulting in more informed policy discovery.

\section{Summary and conclusion}
\label{sec: summary and conclusion}
In this article we have shown that the stochastic optimization approach to importance sampling can be interpreted as a reinforcement learning problem. After presenting the importance sampling problem and a brief introduction to reinforcement learning we have shown how to formulate a MDP for the corresponding stochastic control problem. The MDP is the basic framework for reinforcement learning. By constructing the MDP we have established a first link. We then compared the optimization approaches given for both problems. The comparison has shown that the two optimization approaches are similar and that the optimization in the SOC case is a special case of the reinforcement learning formulation. In the SOC case the forward model of the controlled dynamical system is explicitly given while the reinforcement learning formulation is more general. A third connection has been shown by a detailed discussion of the algorithms developed for the SOC case. Here we have shown that a gradient-based method already proposed in the stochastic optimal control literature can be interpreted as the deterministic policy version of the well-known REINFORCE algorithm which turns out to be model-based. All in all, we have made three connections. We have introduced ideas from reinforcement learning that can be applied to problems seeking optimal deterministic policies namely DPG and its most popular variants DDPG and TD3. These algorithms are model-free policy gradient methods. They can be applied to the importance sampling SOC problem. We have presented the application of both algorithms used in a small dimensional setting and discussed their possible advantages and disadvantages. By applying TD3 to the SOC problem we have clearly shown that the importance sampling SOC problem can be interpreted as a reinforcement learning problem. \\

The main advantage of this is that ideas from reinforcement learning can now be applied to the stochastic optimal control approach to solving the importance sampling problem. For example reinforcement learning has already addressed the question of how to deal with the variance of the gradient estimator. The actor-critic method has been developed to solve this problem and it has been shown that the method achieves this goal. Especially for problems where the time evolution of the dynamical system is strongly influenced by metastable behaviour this can be very helpful to reduce the sampling effort. Furthermore, the issue of efficient data usage has been discussed in the reinforcement learning community and various offline methods have been proposed to solve this problem. 
There are many other interesting ideas that have been addressed by the reinforcement learning community. Thus, this link can be used to efficiently design good and robust algorithms for high-dimensional settings of the importance sampling application. \\

We think that a combination of our model-based gradient estimator with an actor-critic design could be very interesting for the development of algorithms with fast convergence. Another research direction for us is the application of importance sampling to high-dimensional problems such as molecular dynamics. There is already related work exploring these ideas (see, e.g., \cite{Lelievre2022}). However, a stable application to real molecules is still lacking in the literature and would be a very helpful area of application. Another interesting line of research is the combination of model-free and model-based methods. As we have seen in the experiments with higher metastability learning with TD3 become unstable at a certain point. One could use TD3 to compute a good starting point so that the metastability is reduced and then switch to model-based optimization which seems to be much more stable. Similar ideas with pre-initialisation have been proposed in \cite{RiberaQuerRichter2022} where the optimization procedure is combined with an adapted version of the metadynamics algorithm.

\section*{Acknowledgement}
The authors would like to thank the HPC Service of ZEDAT, Freie Universität Berlin, for computing time \cite{Bennett2020}. The research of J.Q. has been funded by the Einstein Foundation Berlin. 
The research of E.R.B has been funded by Deutsche Forschungsgemeinschaft (DFG) through grant CRC 1114 ``Scaling Cascades in Complex Systems'', Project number 235221301, Project A05 ``Probing scales in equilibrated systems by optimal nonequilibrium forcing''. 

\section*{Data availability}
The code used for the numerical examples is available on GitHub at  \url{www.github.com/riberaborrell/rl-sde-is}.

\bibliographystyle{plain}
\bibliography{references}
\end{document}